\documentclass[10pt]{article}

\usepackage[english]{babel}
\usepackage{amssymb}
\usepackage{amsfonts}
\usepackage{cite}
\usepackage{stmaryrd}
\usepackage{mathtools}
\renewcommand{\baselinestretch}{1.3}
\makeatletter
\def\@begintheorem#1#2{\trivlist%
 \item[\hskip \labelsep{\sffamily\bfseries #2\ #1}]\itshape}
\newtheorem{teo}{Theorem}[section]
\newtheorem{defi}[teo]{Definition}

\newtheorem{lem}[teo]{Lemma}
\newtheorem{pro}[teo]{Proposition}
\newtheorem{prop}[teo]{Proposition}
\newtheorem{_rem}[teo]{Remark}
\newtheorem{_rems}[teo]{Remarks}
\newtheorem{_eje}[teo]{Example}

\newenvironment{rem}{\def\@begintheorem##1##2{\trivlist%
 \item[\hskip\labelsep{\sffamily\bfseries ##2\ ##1}]}\begin{_rem}}{\end{_rem}}

\makeatother

\newenvironment{beweis}{{\em Proof.}}{\hfill $\rule{2mm}{2mm}$
\vspace{3mm}\break}

\DeclareMathAlphabet{\Ma}{U}{msa}{m}{n}
\DeclareMathAlphabet{\Mb}{U}{msb}{m}{n}
\DeclareMathAlphabet{\Meuf}{U}{euf}{m}{n}

\def\got#1{\Meuf{#1}}

\DeclareSymbolFont{ASMa}{U}{msa}{m}{n}
\DeclareSymbolFont{ASMb}{U}{msb}{m}{n}
\DeclareMathSymbol{\hrist}{\mathord}{ASMa}{"16}
\DeclareMathSymbol{\varkappa}{\mathalpha}{ASMb}{"7B}
\DeclareMathSymbol{\CrPr}{\mathord}{ASMb}{"6F}

\def\restriction{\upharpoonright}

  \def\al #1.{{\mathcal{#1}}}
  \def\ot #1.{{\got{#1}}}

  \def\ccr #1,#2.{\overline{\Delta(#1,\,#2)}}
  
  \def\wp{\got S}
  \def\b #1.{{\bf #1}}
  \def\cross#1.{\mathrel{\mathop{\times}\limits_{#1}}}
  
  \def\C{\Mb{C}}
  \def\N{\Mb{N}}
  
  \def\R{\Mb{R}}
  \def\Z{\Mb{Z}}
  \def\T{{\Bbb T}}

  \def\rest{\upharpoonright}

\def\f #1,#2.{\mathsurround=0pt \hbox{${#1\over #2}$}\mathsurround=5pt}
\def\bx{{\bf x}}
\def\bg{{\bf g}}
\def\bt{{\bf t}}
\def\by{{\bf y}}
\def\bff{{\bf f}}

  \def\wt{\widetilde}
\def\ilim{\mathop{{\rm lim}}\limits_{\longrightarrow}}
  \def\cross #1.{\mathrel{\raise 3pt\hbox{$\mathop\times\limits_{#1}$}}}
  \def\ol #1.{\overline{#1}}
\def\b #1.{{\bf #1}}
\def\slim{\mathop{\hbox{\rm s-lim}}}
\def\ker{{\rm Ker}\,}

\def\s #1.{_{\smash{\lower2pt\hbox{\mathsurround=0pt $\scriptstyle #1$}}\mathsurround=5pt}}
\def\set #1,#2.{\left\{\,#1\;\bigm|\;#2\,\right\}}
\def\maprightu #1;{\smash{\mathop{\longrightarrow}\limits^{#1}}}
\def\maprightd #1;{\smash{\mathop{\longrightarrow}\limits_{#1}}}
\def\maprightt #1,#2.{\mathrel{\smash{\mathop{\longrightarrow}\limits_{#1}^{#2}}}}
\def\chop{\hfill\break}

\def\XP#1!{\renewcommand{\baselinestretch}{.7}\marginpar{$\leftarrow${\footnotesize #1}\hfil}
 \renewcommand{\baselinestretch}{1}}
\def\XB{\marginpar{
{\footnotesize\bf Change~starts-----}\lower 11pt\hbox{\mathsurround=0pt$
\!\!\displaystyle{
\Bigg\downarrow}$\mathsurround=3pt}}}
\def\XE{\marginpar{{\footnotesize\bf Change~ends-----}\raise 10pt\hbox{\mathsurround=0pt$
\!\!\displaystyle{
\Bigg\downarrow}$\mathsurround=3pt}}}
\def\f #1,#2.{\mathsurround=0pt \hbox{${#1\over #2}$}\mathsurround=5pt}
\def\br#1.{\llbracket #1 \rrbracket}
\def\bbrk #1,#2.{{\langle #1 ,\,#2\rangle}}
\def\eqc#1.{{[#1]\s{\scriptscriptstyle{\mathop{\sim}}}.}}
\def\eqcc#1.{{[#1]\s{\scriptscriptstyle{\mathop{\approx}}}.}}

\def\margin #1.{\marginpar{#1}}

\newcommand{\1}{\mathbf{1}}
\renewcommand{\:}{\colon}
\renewcommand{\hat}{\widehat}
\newcommand{\eset}{\emptyset}
\renewcommand{\phi}{\varphi}
\renewcommand{\tilde}{\widetilde}
\newcommand{\subeq}{\subseteq}
\newcommand{\supeq}{\supseteq}
\newcommand{\la}{\langle}
\newcommand{\ra}{\rangle}
\newcommand{\bra}{\br{\b a.}.}
\newcommand{\brf}{\br{\b f.}.}
\newcommand{\brx}{\br{\b x.}.}
\newcommand{\brxi}{\br{{\bf x}_i}.}
\newcommand{\brxj}{\br{{\bf x}_j}.}
\newcommand{\bry}{\br{\b y.}.}
\newcommand{\im}{\mathop{{\rm im}}\nolimits}
\def\Rep{\mathop{\rm Rep}\nolimits}
\def\into{\hookrightarrow}
\def\Hom{\mathop{\rm Hom}\nolimits}

\def\ba{{\bf a}}
\def\bb{{\bf b}}

\newcommand{\cA}{\mathcal{A}}
\newcommand{\cH}{\mathcal{H}}
\newcommand{\cL}{\mathcal{L}}
\newcommand{\cN}{\mathcal{N}}
\newcommand{\cV}{\mathcal{V}}

\newtheorem{defin}[teo]{Definition}
\newenvironment{defn}{\begin{defin}\rm}{\end{defin}}

\newenvironment{prf}{\noindent{\em Proof:}}{\hfill $\rule{2mm}{2mm}$
\vspace{3mm}\break}

\title{\bf
Infinite Tensor Products of $C_0(\R)$: \\
Towards a Group Algebra for $\R^{(\N)}.$}

\author{
{\sc Hendrik Grundling}  \\[1mm]
{\footnotesize Department of Mathematics,}                           \\
 {\footnotesize University of New South Wales, }    \\
  {\footnotesize Sydney, NSW 2052, Australia.}                    \\
  {\footnotesize hendrik@maths.unsw.edu.au}                         \\
  {\footnotesize FAX: +61-2-93857123} \\
\and
   {\sc Karl--Hermann Neeb}                                      \\[1mm]
 {\footnotesize Fachbereich Mathematik,}                       \\
 {\footnotesize Technische Universit\"at Darmstadt,}                      \\
   {\footnotesize  Schlossgartenstrasse 7,}                            \\
{\footnotesize D--64289 Darmstadt Germany.}                                \\
  {\footnotesize neeb@mathematik.tu-darmstadt.de}}

\date{Running title: Towards a Group Algebra for $\R^{(\N)}.$}

\begin{document}
\maketitle

\begin{abstract}
\noindent
The construction of an infinite tensor product of the C*-algebra $C_0(\R)$ is not
obvious, because it is
nonunital, and it has no nonzero projection. Based on a choice of an approximate identity,
we construct here an infinite tensor product of  $C_0(\R),$
denoted $\al L.\s{\al V.}.,$ and use it to
find (partial) group algebras for the full continuous representation theory of $\R^{(\N)}.$ We obtain
an interpretation of the Bochner--Minlos theorem in $\R^{(\N)}$ as the pure state space decomposition
of the partial group algebras which generate $\al L.\s{\al V.}..$ We analyze the representation theory
of $\al L.\s{\al V.}.,$  and show that there is a bijection between a natural set of representations
of $\al L.\s{\al V.}.$ and
 ${\rm Rep}\big(\R^{(\N)},\al H.\big)\,,$ but that there is an extra part which
 essentially consists of the representation theory of a multiplicative semigroup
 $\al Q.$ which depends on the initial choice of approximate identity.
\end{abstract}


\medskip
\noindent {\bf Mathematics Classification:} 22D25, 46L06, 43A35.

\section{Introduction}

The class of locally compact groups has a rich structure theory with a great many tools developed
to analyze the representation theory of such  groups, e.g.,
 group C*-algebras, induction,
integral decompositions etc. Unfortunately there are many non-locally compact groups which naturally
arise in analysis or physics applications, e.g. mapping groups or inductive limit groups,
and for such groups these tools fail, and one has to do the analysis on a case-by-case basis,
with no systematic theory to draw on.

Here we want to consider the question of how to generalize
the notion of a (twisted) group algebra to topological groups which are not locally compact
(hence have no Haar measure).
Such a generalization,
called a {\it full host algebra,} has been proposed in~\cite{Gr05}. Briefly,
it is a $C^*$-algebra $\al A.$ whose multiplier algebra
$M(\al A.)$ admits a homomorphism $\eta \: G \to U(M({\cal A})),$
such that the (unique) extension of
the representation theory of $\al A.$ to $M(\al A.)$ pulls back
via $\eta$ to the continuous unitary representation theory
of $G$. There is also an analogous concept for
unitary $\sigma\hbox{--representations,}$ where $\sigma$ is a
continuous $\T$-valued $2$-cocycle on $G$.
Thus, given a full host algebra $\al A.,$
the continuous unitary representation theory of $G$ can be analyzed on
$\al A.$ with a large arsenal of $C^*$-algebraic tools.
Such a host algebra need not exist for a general topological group
because there exist topological groups with faithful unitary representations
but without non-trivial irreducible ones (cf.\ \cite{GN01}).
One example of a full host algebra for a group which is not locally compact,
has been constructed explicitly for the
$\sigma\hbox{--representations}$ of an infinite dimensional
topological linear space $S,$ considered as a group  cf.~\cite{GrNe}.

Probably the simplest infinite dimensional group is $\R^{(\N)}$
(the set of real-valued sequences with only finitely many nonzero entries)
with the inductive limit topology w.r.t. the natural inclusions
$\R^n\subset\R^{(\N)}.$ This group is well--studied in stochastic analysis, and will be the
main object of study also in this paper.
Our aim here is to construct explicitly  C*-algebras which have
useful host algebra properties for $\R^{(\N)}.$
Recall that for the group C*-algebras we have:
\[
C^*(\R^n)\otimes C^*(\R^m)\cong C^*(\R^{n+m})
\]
and this suggests that for a host algebra of $\R^{(\N)}$ we should try an
infinite tensor product of $C^*(\R).$ This is difficult to do, for two reasons:
\begin{itemize}
\item{} $C^*(\R)\cong C_0(\R)$ is nonunital, and the standard infinite tensor products
of C*-algebras require unital algebras.
\item{} There is a definition for an infinite tensor product of nonunital algebras developed by
Blackadar cf.~\cite{Bla}, but this requires the algebras to have nonzero projections, and the
construction depends on the choice of projections.
(We used this construction in \cite{GrNe} to construct an infinite tensor product to produce a host algebra.)
However, $C^*(\R)\cong C_0(\R)$ has no nonzero projections, so this method will not work.
\end{itemize}
\noindent In the light of these difficulties, we will develop here an
infinite tensor product of $C_0(\R)$ relative to a choice of  approximate identity in each entry, to
replace the choice of projections in Blackadar's approach. As expected, the construction
will depend on the choice of approximate identities, though it still produces for each choice an
algebra with strong host algebra properties.

The construction of (``semi-'')host algebras for $\R^{(\N)}$ will aid our understanding of the
Bochner--Minlos theorem. We first recall:
\begin{teo}
(Bochner--Minlos Theorem for $\R^{(\N)}$)
There is a bijection between continuous normalized positive definite
functions (states) $\omega$ of
$\R^{(\N)}$ and regular Borel probability measures $\mu$ on $\R^{\N}$ (with product topology) given by the
Fourier transform:
\[
\omega(\b x.)=\int_{\R^{\N}}e^{i\b x.\cdot\b y.}d\mu(\b y.)\,,\quad\b x.\in \R^{(\N)}
\]
where $\b x.\cdot\b y.:=\sum\limits_{n=1}^\infty x_ny_n\,,$ $\b x.\in \R^{(\N)},$ $\b y.\in \R^{\N}.$
\end{teo}
If we replace both $\R^{(\N)}$ and $\R^{\N}$ by $\R^n,$ this is the classical Bochner theorem, which we can obtain
immediately from the state space integral decomposition of any state of $C^*(\R^n)\cong C_0(\R^n)$ in terms of pure states.
This suggests that if we have a host algebra of $\R^{(\N)},$ we can obtain the Bochner--Minlos theorem
from state space decompositions of states on the host algebra in terms of pure states.
We will see below that we can already obtain the Bochner--Minlos theorem from the
weaker ``semi--host'' algebras which we will construct.

The structure of this paper is as follows. In Section~\ref{DefsNotn}
we collect the basic definitions and notation for host algebras, in
Section~\ref{AlgTens} we give a detailed treatment of the aspects
of infinite tensor products which we will need for this paper.
 In Section~\ref{GausH} we start in a concrete setting  on $L^2(\R^{\N},\mu), $ where $\mu$
 is a product measure of
 probability measures, each absolutely continuous w.r.t. the Lebesgue measure,
 and we construct an infinite
 tensor product of $C_0(\R)$ w.r.t. a choice (compatible with $\mu)$
 of approximate identity in each entry. This concrete C*-algebra can already
 produce Bochner--Minlos decompositions for the limited class of
 positive definite functions on $\R^{(\N)}$ associated with it.
In Section~\ref{ParHRN} we develop abstractly the infinite
 tensor product of $C_0(\R)$ w.r.t. an arbitrary choice of elements of a
 fixed approximate identity, we analyze its representation theory
 and through the unitary embedding of $\R^{(\N)}$ in its multiplier algebra,
 we consider the relation of its representation theory to that of
 $\R^{(\N)}.$ We find that it can adequately model a subset of the
 representation theory  of
 $\R^{(\N)},$ but there is a small additional part.
 We show that the Bochner--Minlos decompositions for {\it any}  continuous
 positive definite function on $\R^{(\N)}$ can be obtained from the pure state space
 decomposition of these algebras. Finally, in Section~\ref{FullRepRN}, we collect
 these algebras together in one large C*-algebra, which we show, can model the full
 continuous representation theory  of $\R^{(\N)}.$
 However, the representation theory of this algebra also has an additional part  which
 essentially consists of the representation theory of a multiplicative semigroup
 $\al Q.$ which depends on the initial fixed choice of approximate identity.

\section{Definitions and notation}
\label{DefsNotn}

We will need the following notation and concepts for our main results.
\begin{itemize}
\item{}
In the following, we write $M({\cal A})$ for the
multiplier algebra of a $C^*$-algebra ${\cal A}$ and, if ${\cal A}$ has a unit, $U({\cal A})$ for its
unitary group.
We have an injective morphism of $C^*$-algebras
$\iota_{\cal A} \: {\cal A} \to M({\cal A})$ and will just denote ${\cal A}$ for its
image in $M({\cal A})$. Then ${\cal A}$ is dense
in  $M({\cal A})$ with respect to
the {\it strict topology},
which is the locally convex topology defined by the seminorms
$$ p_a(m) := \|m \cdot a\| + \|a \cdot m\|,
\qquad a\in {\cal A},\; m\in M({\cal A})$$
(cf.\cite{Wor1}).
\item{}
For a complex Hilbert space ${\cal H}$, we write $\Rep({\cal A},{\cal H})$ for the
set of non-degenerate representations of ${\cal A}$ on ${\cal H}$.
Note that the collection $\Rep {\cal A}$ of all non-degenerate
representations of ${\cal A}$ is not a set, but a (proper) class
in the sense of von Neumann--Bernays--G\"odel set theory, cf.~\cite{TZ75}, and in
this framework we can consistently manipulate the object $\Rep {\cal A}.$
However, to avoid set--theoretical subtleties, we will express our results below
concretely, i.e., in terms of $\Rep({\cal A},{\cal H})$ for given Hilbert spaces
$\al H..$
We have an injection
$$ \Rep({\cal A}, {\cal H}) \into \Rep(M({\cal A}),{\cal H}), \quad \pi \mapsto \tilde\pi
\quad \hbox{ with } \quad \tilde\pi \circ \iota_{\cal A} = \pi, $$
which identifies the non-degenerate representation $\pi$ of
${\cal A}$ with that representation $\tilde\pi$ of its multiplier algebra which
extends $\pi$ and is
continuous with respect to the  strict topology on $M({\cal A})$ and
the topology of pointwise convergence on $B({\cal H})$.
We will refer to $\tilde\pi$ as the {\it strict extension} of $\pi,$ and it is
easily obtained by
\[
\wt{\pi}(M)=\slim_{\lambda\to\infty}\pi(M E_\lambda)\quad\forall\,M\in M(\al A.)
\]
where $\{E_\lambda\}_{\lambda\in\Lambda}\subset\al A.$ is any approximate
identity  of $\al A..$
\item{}
For  topological groups $G$ and $H$ we write
$\Hom(G,H)$ for the set of continuous group homomorphisms $G \to H$.
We also write
$\Rep(G,{\cal H})$ for the set of all (strong operator)
continuous unitary representations of
$G$ on ${\cal H}\,.$ 
Endowing $U({\cal H})$ with the strong operator topology
turns it into a topological group, denoted $U({\cal H})_s$,
so that $\Rep(G,{\cal H}) = \Hom(G,U({\cal H})_s)$.
The set of continuous normalized positive definite functions on $G$ (also
called {\it states}) and denoted by $\wp(G),$
is in bijection with the state space of the
group C*-algebra $C^*(G)$  when $G$ is locally compact.
If $G$ is not locally compact, $\wp(G)$ is in bijection with a
subset of the state space of  $C^*(G_d)$, where $G_d$ denotes $G$ with the discrete topology,
and the question arises as to whether there is a C*-algebra which can play the role of
$C^*(G).$ We clarify first what is meant by this:
\end{itemize}
\begin{defn}
Let  $G$ be a topological group.\chop
A {\it host algebra for  $G$} is a pair
$({\cal L}, \eta)$ where ${\cal L}$ is a $C^*$-algebra
and $\eta \: G \to U(M({\cal L}))$ is a homomorphism
such that for each complex Hilbert space
${\cal H}$ the corresponding map
$$ \eta^* \: \Rep({\cal L},{\cal H}) \to \Rep(G, {\cal H}), \quad
\pi \mapsto \tilde\pi \circ \eta $$
is injective. We then write $\Rep(G,{\cal H})_\eta \subeq \Rep(G,{\cal H})$ for the range
of $\eta^*$.
We say that $(\cL,\eta)$ is a {\it full host algebra} of $G$ if
$\eta^*$ is surjective for each Hilbert space~${\cal H}$.
If the map $\eta^*$ is not injective, we will call the pair
$({\cal L}, \eta)$ a {\it semi-host algebra} for $G.$

Note that by the universal property of  group algebras, the homomorphism
$\eta \: G\to U(M({\cal L}))$  extends uniquely to the
discrete group C*-algebra $C^*(G_d),$ i.e. we have a *-homomorphism
$\eta \: C^*(G_d)\to U(M({\cal L}))$ (still denoted by $\eta$).

A similar notion can also be defined for projective representations
(cf.~\cite{GrNe}).
\end{defn}

\begin{rem}
\begin{itemize}
\item[(1)]
It is well known that for each locally compact group $G$,
the group $C^*$-algebra $C^*(G)$, and
the natural map $\eta_G \: G \to M(C^*(G))$ provide a
full host algebra (\cite[Sect.~13.9]{Dix}). The map $\eta_G: G\to M(C^*(G))$
is continuous w.r.t. the strict topology of
$M(C^*(G))$ (this is an easy consequence of the fact that
$\im(\eta_G)$ is bounded and that the action on the corresponding $L^1$-algebra
is continuous).
\item[(2)]
Note that for a host algebra $(\cL,\eta)$
the map $\eta^*$ preserves direct sums, unitary conjugation, subrepresentations,
and for full host algebras, irreducibility (cf.~\cite{Gr05}).
\item[(3)]
When $({\cal L}, \eta)$ is merely a semi-host algebra  for $G,$ then the map $\eta^*$ still
preserves direct sums, unitary conjugation, subrepresentations, but in general, not irreducibility.
However, in the case that $G$ is Abelian (as it will be in this paper), since irreducible representations are
just characters, and the map $\eta^*$ takes one--dimensional representations to one--dimensional ones,
here it will preserve irreducibility. So for Abelian groups, semi--hosts are useful to carry
representation structure (e.g. integral decompositions) from the representation theory of $\al L.$ to
the representation theory of $G,$ and we will use that in this paper to analyze the Bochner--Minlos theorem.
\end{itemize}
\end{rem}

\section{Basic Theory of Infinite Tensor Products}
\label{AlgTens}

Since we need to develop the concept of infinite tensor products of non-unital algebras,
it is necessary to collect first some basic material on infinite tensor products,
and to fix notation.
We follow Bourbaki~\cite{Bou89} and Wegge--Olsen~\cite{WO}.
There are several different concepts of infinite tensor products of
unital algebras. See Bourbaki~\cite{Bou89}, Guichardet~\cite{Gui}, Araki~\cite{Ara},
though infinite tensor products of algebras without identity are only
done in Blackadar~\cite{Bla}.

\subsection{Algebraic tensor products of arbitrary many factors.}

\begin{defn} Let $(X_t)_{t \in T}$ be an indexed set of non-zero complex vector
spaces,
where $T$ can have any cardinality.
We write $\bx = (x_t)_{t \in T}$ for the elements of the product
space $\prod\limits_{t \in T}X_t.$
A map $f \: \prod\limits_{t \in T} X_t \to V$ to a vector space
$V$ is said to be {\it multilinear} if it is linear in each entry. That is, for each
$t_0\in T$ and $\bx\in\prod\limits_{t \in T\backslash \{t_0\}}X_t,$
the map
$$ X_{t_0} \to V, \quad y_{t_0} \mapsto f(\bx\times y_{t_0}) $$
is linear, where $\bx\times y_{t_0}=:{\bf z}\in\prod\limits_{t \in T}X_t$
is that element for which $z_t = x_t$ if $t\not= t_0$ and
$z_{t_0}=y_{t_0}.$

A pair $(\iota, V)$ consisting of a vector space $V$ and a multilinear map
$\iota \: \prod\limits_{t \in T} X_t \to V$ is called an  {\it (algebraic)
tensor product
of  $(X_t)_{t \in T}$} if it has the following universal
property
\begin{description}
\item[\bf(UP)] For each multilinear map $\phi \: \prod\limits_{t \in T} X_t \to W$, there exists a unique linear map $\tilde\phi \: V \to W$ with
$\tilde\phi \circ \iota = \phi$.
\end{description}

The usual arguments (cf. Proposition~T.2.1 \cite{WO}) show that the universal property determines
a tensor product up to linear isomorphism (factoring through the maps $\iota).$
 We may thus denote $V$ by
$\bigotimes\limits_{t \in T} X_t$
and denote the {\it elementary tensors} by
$$ \mathop{\otimes}_{t \in T} x_t := \iota(\bx)\in \bigotimes_{t \in T} X_t,\quad
\hbox{for}\quad \bx\in
\prod_{t \in T} X_t\,.$$
\end{defn}
To simplify notation, we write $X := \prod\limits_{t \in T} X_t$ in the following.
Observe that no order in $T$ appears in this definition, so e.g.
$X_1\otimes X_2$ and $X_2\otimes X_1$ (in the usual notation) will be identified.
\begin{lem} For each indexed set $(X_t)_{t \in T}$ of complex vector spaces,
a tensor product $(\iota, \bigotimes\limits_{t \in T} X_t)$ exists.
\end{lem}

\begin{prf} (cf.\ \cite[Ch.~II,\S 3.9]{Bou89} for a more general construction)
We consider the free complex vector space
\[
\C^{(X)}:=\big\{f:X\to\C\,\mid\,{\rm supp}(f)\;\hbox{is finite}\big\}
={\rm Span}\big\{\delta_{\bf x}\,\mid\,{\bf x}\in X\big\}
\]
where  $\delta_{\bf x}({\bf y})=
1$ if ${\bf x} = {\bf y}$ and zero otherwise.
Note that $\big\{\delta_{\bf x}\,\mid\,{\bf x}\in X\big\}$ is a basis
for $\C^{(X)}$. Define the sets
\begin{eqnarray*}
N_a &:=&\big\{ \delta_{\bf x}+\delta_{\bf y}-\delta_{\bf z}\,\big|\,
\exists\,r\in T\;\hbox{such that}\; x_r+y_r=z_r\,,\;\hbox{and} 
\; x_t=y_t=z_t\;\;\forall\,t\not=r\big\} \\[1mm]
N_m &:=&\big\{ \delta_{\bf x}-\mu\delta_{\bf y}\,\big|\,\mu\in\C,\;\hbox{and}\;
\exists\,r\in T\;\hbox{such that}\; x_r=\mu y_r\,,\;\hbox{and}   
\; x_t=y_t\;\;\forall\,t\not=r\big\}  \\[1mm]
\al N.&:=& {\rm Span}\big(N_a\cup N_m\big)\subset\C^{(X)}\;.
\end{eqnarray*}

We now consider the quotient space
$V := \C^{(X)}\big/\cN$ and write
$\iota \: X \to V, \bx \mapsto \delta_\bx + \cN$
for the induced map. The definition of $\cN$ immediately implies
that $\iota$ is multilinear and we only have to verify the universal
property.

Let $\phi \: X \to M$ be a multilinear map.
We extend $\phi$ to a linear map $\phi:\C^{(X)}\to M$ by
$\varphi(f):=\sum\limits_{\b x.\in X}f(\b x.)\,\varphi(\b x.)$.
The multilinearity of $\phi$ now implies that its linear extension
annihilates the subspace $\cN,$ hence it factors through a linear map
$\tilde\phi \: V \to M$ satisfying
$\tilde\phi \circ \iota = \phi.$ That $\tilde\phi$ is uniquely determined
by this property follows from the fact that
$\im(\iota)$ spans $V.$
\end{prf}
\begin{teo}
\label{Assoc}
(Associativity) \chop
Let $\set T_s\subset T,s\in S.$ be a partition of $T$ such that $|S|<\infty\,.$ Then the map
$$ \psi : \prod_{t \in T} X_t \to
\mathop{\textstyle\bigotimes}_{s\in S}\big(\mathop{\textstyle\bigotimes}_{t_s\in T_s}X_{t_s}\big),
\quad \psi((x_t)_{t \in T})
:= \mathop{\textstyle\bigotimes}_{s\in S}
\big(\mathop{\textstyle\bigotimes}_{t_s\in T_s}x_{t_s}\big) $$
is multilinear and factors through a linear isomorphism
$\tilde\psi \: \bigotimes\limits_{t \in T} X_t \to
\mathop{\textstyle\bigotimes}\limits_{s\in S}\big(\mathop{\textstyle\bigotimes}\limits_{t_s\in T_s}X_{t_s}\big).$
\end{teo}
\begin{beweis} It is clear from the definition that $\psi$ is multilinear,
so we obtain a unique linear map $\wt{\psi}:\mathop{\bigotimes}\limits_{t\in T}X_t
\to\mathop{\bigotimes}\limits_{s\in S}\big(\mathop{\bigotimes}\limits_{t\in T_s}X_t\big)$
with $\tilde\psi \circ \iota = \psi$.

To see that $\tilde\psi$ is a linear isomorphism, it suffices to observe that
the multilinear map $\psi$ has the universal property (UP).
So let $\phi \: X \to V$ be a multilinear map.
With $Y_s := \prod\limits_{t \in T_s} X_t$, we have
$X = \prod\limits_{s \in S} Y_s$. Then for each
$s_0\in S$ and for each
$\b y.\in\prod\limits_{s \in S\backslash s_0} Y_s$ we obtain
a unique  map
$$ \varphi^{s_0}_{\b y.}:Y_{s_0}=\prod\limits_{t \in T_{s_0}} X_t \to V, \quad \varphi^{s_0}_{\b y.}(y_{s_0}):=  \phi(\b y.\times y_{s_0})\,,$$
which is clearly multilinear w.r.t. the factors $\prod\limits_{t \in T_{s_0}} X_t=Y_{s_0}$ hence induces a linear map on
$\bigotimes\limits_{t \in T_{s_0}} X_t.$
Since $\b y.\mapsto\varphi^{s_0}_{\b y.}(v)$ is  multilinear in $\b y.\in\prod\limits_{s \in S\backslash s_0} Y_s$
for fixed $v\in\bigotimes\limits_{t \in T_{s_0}} X_t,$
we can apply the argument again to an  $s_1\not= s_0\in S$ for this map, and then
continue the process until we have exhausted $S.$
This produces a multilinear map
$$\hat \phi \: \prod_{s \in S} \Big(\bigotimes_{t \in T_s} X_t\Big) \to V $$
which factors through a linear map
$$\tilde \phi \: \bigotimes_{s \in S} \Big(\bigotimes_{t \in T_s} X_t\Big) \to V
\quad \mbox{ with } \quad
\tilde\phi\big(\mathop{\otimes}_{s \in S}\big(\mathop{\otimes}_{t_s \in T_s} x_{t_s}\big)\big)
= \phi((x_t)_{t \in T}),$$
i.e., $\tilde\phi \circ \psi = \phi$.
Moreover, since $\bigotimes\limits_{s \in S} \Big(\bigotimes\limits_{t \in T_s} X_t\Big)$ is spanned
by elements of the form $\mathop{\otimes}\limits_{s \in S}\big(\mathop{\otimes}\limits_{t_s \in T_s} x_{t_s}\big)$ it follows that
$\tilde \phi$ is uniquely determined by the last equation.
Thus $\psi$ has the universal property (UP), hence
$\tilde\psi$ is a linear isomorphism.
\end{beweis}

\begin{rem}
Associativity does not seem to hold for a partition of $T$ into infinitely many sets (i.e., for
$|S|=\infty$). This is because $\mathop{\bigotimes}\limits_{t\in T}X_t$ is spanned
by elementary tensors, and ${\mathop{\bigotimes}\limits_{s\in S}
\Big(\sum\limits_{t_s=1}^{n_s}\mathop{\bigotimes}\limits_{r_s\in T_s}x_{r_s}^{(t_s)}\Big)}$ cannot
be written as a finite linear combination of elementary tensors if there are infinitely many $s\in
S$ with $n_s>1\,.$
\end{rem}

\begin{defn} (a) Assume that $(X_t)_{t \in T}$ is a family of complex algebras.
We now construct an algebra structure on their tensor product.
For each fixed $\b x.\in X=\prod\limits_{t\in T}X_t$, define a map
\[
\mu_{\b x.}:X\to\mathop{\textstyle\bigotimes}_{t\in T}X_t\quad\hbox{by}\quad
\mu_{\b x.}(\b y.):=\mathop{\textstyle\bigotimes}_{t\in T}x_ty_t =\iota(\b x.\cdot\b y.)
\]
where $\b x.\cdot\b y.\in X$ is given by $(\b x.\cdot\b y.)_t:=x_ty_t$ for all $t\in T\,,$
and we will also let $\b x.^n\in X$ denote $(\b x.^n)_t:=(x_t)^n$ for all $t\in T$ and $n\in\N.$
Since $\mu_{\b x.}$ is multilinear, it induces a linear map
$$\mu_{\b x.}:\mathop{\bigotimes}\limits_{t\in T}X_t\to\mathop{\bigotimes}\limits_{t\in T}X_t.$$
This defines a multilinear map
\[
\mu:X\to{\rm End}\big(\mathop{\textstyle\bigotimes}_{t\in T}X_t\big)\quad\hbox{by}\quad
\mu(\b x.):=\mu_{\b x.}
\]
and thus a linear map
$\mu:\mathop{\bigotimes}\limits_{t\in T}X_t\to{\rm End}\big(\mathop{\bigotimes}\limits_{t\in T}X_t\big)\,.$
Explicitly we have for $a=\sum\limits_i\iota(\b x._i)$ and $b=\sum\limits_j\iota(\b y._j)\in
\mathop{\bigotimes}\limits_{t\in T}X_t$ that
\[
\mu(a)(b)=\sum_i\mu_{\b x._i}\Big(\sum_j\iota(\b y._j)\Big)
=\sum_i\sum_j\mu_{\b x._i}\big(\iota(\b y._j)\big)=\sum_i\sum_j\iota({\b x._i}\cdot\b y._j)
\]
where the sums are finite. We denote the multiplication as usual by
$a\,b:=\mu(a)(b)$ for $a,\,b\in\mathop{\bigotimes}\limits_{t\in T}X_t$.
Associativity for this multiplication follows from componentwise associativity,
and hence $\mathop{\bigotimes}\limits_{t\in T}X_t$ is an algebra over $\C\,.$

(b) Next, we assume, in addition, that each $X_t$ is a $*$-algebra.
We want to turn $\mathop{\bigotimes}\limits_{t\in T}X_t$ into a $*$-algebra.
Given any vector space $V$ over $\C\,,$ let $V^c$ denote the conjugate vector space.
Thus, for each $t\in T$, the involution $*:X_t\to X_t^c$ becomes a $\C\hbox{--linear}$
map (instead of conjugate linear on $X_t$).
Define a map
\[
\gamma:X\to\big(\mathop{\textstyle\bigotimes}_{t\in T}X_t\big)^c\quad\hbox{by}\quad
\gamma(\b x.):=\mathop{\textstyle\bigotimes}_{t\in T}x_t^*=\iota(\b x.^*)
\]
where $\b x.^*\in X$ is given by $(\b x.^*)_t:=x^*_t$ for all $t\in T\,.$
Since $\gamma$ is multilinear, it defines a
linear map
$\gamma:\mathop{\bigotimes}\limits_{t\in T}X_t\to\big(\mathop{\bigotimes}\limits_{t\in T}X_t\big)^c$.
Its intertwining properties with multiplication then follow from the
componentwise properties. As usual,
we write $a^* := \gamma(a)$ for $a\in\mathop{\bigotimes}\limits_{t\in T}X_t\,,$
and hence $\mathop{\bigotimes}\limits_{t\in T}X_t$ becomes
a $*$-algebra over $\C\,.$
\end{defn}

This defines the basic objects which we will work with.

\subsection{Stabilized spaces.}

 We will also need the following structures.
\begin{defn}
\label{Eqce} We define an equivalence relation on $X$ by
$\b x.\sim\b y.$ whenever the set
${\{ t \in T \mid x_t \not=y_t\}}$ is finite.
Denote the equivalence class of $\b x.\in X$ by $\eqc{\b x.}.$ and define
\[
 \br{\b x.}.:={\rm Span}\big\{\mathop{\otimes}_{t\in T}y_t\,\mid\,
\b y.\sim\b x.\big\}\subset\mathop{\textstyle\bigotimes}_{t\in T}X_t\;.
\]
\end{defn}

\begin{prop} \label{prop:2.7} The following assertions hold:
  \begin{description}
  \item[\rm(i)] For any pair  $(\bx,F)$ such that $\bx \in X$ and  $F \subeq T$ a finite subset
with $x_t \not=0$ for $t \not\in F$, there exists a linear map
$$ \phi_F \: \bigotimes_{t \in T} X_t \to \bigotimes_{t \in F} X_t $$
satisfying $\bry \subeq \ker \phi_F$ for $\by \not\sim \bx$ and
$$ \phi_F\big((\mathop{\otimes}_{t \in F} y_t) \otimes (\mathop{\otimes}_{t \not\in F} x_t)\big)
=  \mathop{\otimes}_{t \in F} y_t \quad \mbox{ for } \quad y_t \in X_t, t \in F. $$
  \item[\rm(ii)] $\brx\not= \{0\}$ if and only if at most finitely many
components of $\bx$ vanish.
  \item[\rm(iii)]
The subspace $\brx$ is isomorphic to the direct limit of the
finite tensor products $\mathop{\otimes}\limits_{t \in J} X_t$, $J \subeq T$ finite,
with respect to the
connecting maps
\[
 \varphi\s{K,J}.:\mathop{\textstyle\bigotimes}_{t\in J}X_t\to \mathop{\textstyle\bigotimes}_{t\in K}X_t
\quad\hbox{with}\quad \varphi\s{K,J}.\big(\mathop{\otimes}_{t\in J}y_t\big):=
\big(\mathop{\otimes}_{t\in J}y_t\big)\otimes\big(\mathop{\otimes}_{s\in K\backslash J}x_s\big).
\]
  \item[\rm(iv)]  $\bigotimes\limits_{t \in T} X_t$ is the direct
sum of the subspaces
$\brx$, $\bx \in X$.
  \end{description}
\end{prop}

\begin{prf} (i) For $t \not\in F$ we pick linear functionals
$\lambda_t \in X_t^*$ with $\lambda_t(x_t) = 1$ and define a map
$$ \hat\phi_F \: X \to \bigotimes_{f \in F} X_f, \quad
\hat\phi_F({\bf y}) :=
\begin{cases}
\prod\limits_{t \in T\setminus F} \lambda_t(y_t) \cdot \big(\mathop{\otimes}\limits_{s \in F} y_s\big)
     & \text{for } {\bf y} \sim {\bf x} \\
0  & \text{for } {\bf y}\not\sim {\bf x}.
\end{cases}
$$
We claim that $\hat\phi_F$ is multilinear. To see that $\hat\phi_F$
is linear in the $t$-component, let
$\by, \by' \in X$ with $y_s = y_s'$ for $s\not=t$.
Then either both are equivalent to $\bx$ or none is. In
either case, the definition of $\hat\phi_F$ implies the linearity of the map
$z_t \mapsto \hat\phi_F(\by\times z_t).$
Therefore $\hat\phi_F$ is multilinear, hence induces a linear map
$$  \phi_F \: \bigotimes X_t \to \bigotimes_{t \in F} X_t  $$
satisfying all requirements.

(ii) If the set $\{ t \in T \mid x_t = 0\}$ is finite, then (i)
implies that $\br{\b x.}.\not=\{0\}$ since
none of the spaces $X_t$ vanishes by our initial assumption.
We also note that, if infinitely many $x_t$ vanish,
then $\brx$ is spanned by elements $\iota(\by)$, where
$\by$ has at least one zero entry. Then $\iota(\by) = 0$, and
consequently $\brx = \{0\}$.

(iii) Let $J\subset K\subset T$ such that $|K|<\infty$.
Then we obtain linear maps
\[
 \varphi\s{K,J}.:\mathop{\textstyle\bigotimes}_{t\in J}X_t\to \mathop{\textstyle\bigotimes}_{t\in K}X_t
\quad\hbox{with}\quad \varphi\s{K,J}.\big(\mathop{\otimes}_{t\in J}y_t\big):=
\big(\mathop{\otimes}_{t\in J}y_t\big)\otimes\big(\mathop{\otimes}_{s\in K\backslash J}x_s\big).
\]
Since $\varphi\s{L,K}.\circ\varphi\s{K,J}.=\varphi\s{L,J}.$
for $J\subset K\subset L\,,$ and $|L|<\infty\,,$ this is an inductive system. We write
$\ilim{\big(\mathop{\bigotimes}\limits_{t\in J}X_t,\,\varphi\s{K,J}.\big)}$
for its limit. We also have linear maps
\[
 \varphi\s{J}.:\mathop{\textstyle\bigotimes}_{t\in J}X_t\to \br{\b x.}.
\quad\hbox{by}\quad \varphi\s{J}.\big(\mathop{\otimes}_{t\in J}y_t\big):=
\big(\mathop{\otimes}_{t\in J}y_t\big)\otimes\big(\mathop{\otimes}_{s\in T\backslash J}x_s\big)
\in\brx
\]
satisfying $\phi_K \circ \phi_{K,J} = \phi_J$, so that they induce a
linear map
$\varphi:\ilim{\big(\mathop{\bigotimes}\limits_{t\in J}X_t,\,\varphi\s{K,J}.\big)}\to
\brx$. As every element of $\br{\b x.}.$ lies in the image of some
map $\varphi\s{J}.\,,$ and by (i) this map is injective if
$J \supeq \{ t \in T \mid x_t =0\}$,
$\phi$ is a linear isomorphism.

(iv) Since $\iota(\bx)$ is contained $\brx$, it suffices
to show that the sum of the non-zero subspaces $\brx$ is direct.
Suppose that the elements
$\bx_1,\ldots, \bx_n$ are pairwise non-equivalent with
\break $\brxi\not=\{0\}$,  and
that $v_i \in \brxi$ satisfy $\sum_i v_i = 0$.
From (i) we know that there exists for each
$i$ and each finite subset $F \supeq \{ t \in T  \mid x_{i,t} = 0\}$
a linear map
$$ \phi^{(i)}_{F} \: \bigotimes_{t \in T} X_t \to \bigotimes_{t \in F} X_t
\quad \mbox{ with } \quad
\phi^{(i)}_{F}\big((\mathop{\otimes}_{t \in F} y_t) \otimes
(\mathop{\otimes}_{t \not\in F} x_{i,t})\big)
=  \mathop{\otimes}_{t \in F} y_t $$
and vanishing on $\brxj$ for $j \not=i$.
We conclude that
$\phi^{(i)}_{F}(v_i) =0$ for each $F$. Since $F$ can be chosen
arbitrarily large, the definition of $\brxi$ now implies that $v_i = 0$.
\end{prf}

\begin{rem}
If each $X_t$ is an algebra and
$x_t^2 =x_t$ holds for all but finitely many $t \in T$, then
the linear space $\br{\b x.}.$ is a subalgebra.
If each $X_t$ is a $*$-algebra and
$x_t^*=x_t=x_t^2$ for all but finitely many $t\in T$,
then $\br{\b x.}.$ is a $*$-subalgebra.
In the literature (on topological tensor products),  suitable closures of
$\br{\b x.}.$ are often called stabilized infinite tensor products
(stabilized by $\bx$).
\end{rem}

\begin{rem}
In particular, for
 $\bx, \by \in X$ with $\br{\b x.}.\not=\{0\}\not=\br{\b y.}.\,,$ we
have that $\br{\b x.}.\cap\br{\b y.}. =\{0\}$ if and only if
$\b x.\not\sim\b y.\,.$
So, if
$y_t=\lambda_tx_t$ where $\lambda_t\not=1$ for infinitely many $t\in T\,,$ then $\b x.\not\sim\b y.$ and hence
$\mathop{\otimes}\limits_{t\in T}\lambda_tx_t$ is not a multiple of
$\mathop{\otimes}\limits_{t\in T}x_t.$
This is different in Guichardet's
version~\cite{Gui} of continuous  tensor products.
\end{rem}

When the $X_t$ are algebras, we have the following algebraic relations for the spaces $\br{\b x.}.$ in the algebra
$\mathop{\bigotimes}\limits_{t\in T}X_t\,.$
\begin{teo}
\label{ProdCls} If each $X_t$ is a complex algebra, then
\begin{itemize}
\item[\rm(i)] $\br{\b x.}.\cdot\br{\b y.}.\subseteq\br{\b x.\cdot\b y.}.$ for all $\b x.,\,\b y.\in X\,.$
If $X_t\cdot X_t=X_t$ for all $t\,,$  then we have the equality:
${\rm Span}\big(\br{\b x.}.\cdot\br{\b y.}.\big)=\br{\b x.\cdot\b y.}.\,.$
\item[\rm(ii)] $\br{\b x.}.^*=\br{\b x.^*}.$ for all $\b x.\in X$ if all $X_t$ are $*$-algebras.
\item[\rm(iii)] If $\eset\not=
G_t\subset X_t\setminus \{0\}$
is a nonzero multiplicative semigroup for each $t\in T\,,$ then
\[
\al M.:= \sum_{\ba \in \prod_{t \in T} G_t} \bra  \qquad\hbox{(finite sums)}
\]
is a subalgebra of $\mathop{\bigotimes}\limits_{t\in T}X_t\,.$
If in addition, each $X_t$ is a $*$-algebra and
each $G_t$ is $*$-invariant, then $\al M.$ is a
 $*$-subalgebra.
\end{itemize}
\end{teo}
\begin{beweis}
(i) Since $\brx$ is spanned by elements of the form $\iota(\ba)$,
$\ba \sim \bx$ and $\bry$ likewise by elements $\iota(\bb)$ with
$\bb \sim \by$, and we have $\ba\cdot\bb\sim\bx\cdot\by,$
the first assertion follows from
$\iota(\ba) \iota(\bb) = \iota(\ba \cdot\bb) \in \br{\b x.\cdot\b y.}.\,.$

To show that we have equality when  $X_t\cdot X_t=X_t$
for all $t\,,$ note that
$\br{\b x.\cdot\b y.}.$ is spanned by elements of the
form
$\iota(\b a.)=\Big(\mathop{\otimes}\limits_{s\in S}a_s\Big)
\otimes\Big(\mathop{\otimes}\limits_{t\in T\backslash S}x_ty_t\Big)\,,$ where $S$ is finite.
 Since each $a_s \in X_sX_s$ by assumption, it follows that
$\iota(\ba) \in \brx \bry$,
which proves the required equality.

(ii) Since $*$ is involutive, it suffices to show that
$\brx^* \subeq \br{\b x^*.}.$. As $\brx^*$ is spanned by elements
of the form
$\iota(\ba)^*$, $\ba \sim \bx$, the assertion follows from
$\iota(\ba)^* = \iota(\ba^*)$ with $\ba^* \sim \bx^*$.

(iii) Since the set $\{\b x.\in X\,\mid \,x_t\in G_t\;\forall \,t \in T\}$ is a semigroup  w.r.t.\ the componentwise multiplication,
the first statement regarding $\al M.$ follows from (i).
The second statement likewise follows from (ii).
\end{beweis}

\begin{rem} \label{rem:2.11}
(a) Regarding the condition $X_t\cdot X_t=X_t$ in part~(i), this is easily fulfilled, since by
Theorem~5.2.2 in~\cite{Pa94}, we know that
if $\al A.$ is a Banach algebra with a bounded
left approximate identity and $T : \al A.\to \al B.(X)$ is a continuous representation
of $\al A.$ on the Banach space $X\,,$ then for each $y\in\overline{{\rm Span}(T(\al A.)X)}$ there are
elements $a\in\al A.$ and $x \in X$ with $y = T(a)x\,.$
Thus, if $X =\al A.$ and $T : \al A.\to \al B.(X)$ is defined by
$T(A)B:=AB$, then since $\al A.$
has an approximate identity, we have $\al A.=\overline{{\rm Span}(T(\al A.)X)}$ and hence
$\al A.\cdot\al A.=\al A.\,.$ In particular, $\al A.\cdot\al A.=\al A.$ for any $C^*$-algebra $\al A.\,.$

(b)
In regard to the choice of semigroup $G_t$ in (iii) above, when one has unital algebras, the conventional
choice is to set all $G_t=\{\1\}$.
If the $*$-algebras $X_t$ are nonunital but have projections, then
one can take each $G_t$ to be a projection (cf. Blackadar~\cite{Bla}) though the final tensor product algebra
depends on this choice of projections. If the $*$-algebras $X_t$ have no nonzero projections, e.g.~$C_0(\R)$ below,
then we will choose each $G_t$ to be a small $*$-closed semigroup generated by one element (which will be positive, of norm $1$).
\end{rem}

\subsection{Tensor products of representations.}

Below, we will need to complete some $*$-subalgebras of the algebraic tensor product in the operator norm
of a suitable representation, hence need to make explicit the structures involved with infinite tensor
products of Hilbert space representations.

Let $(\cH_t)_{t \in T}$ be a family of Hilbert spaces.
We want to equip selected subspaces of $\mathop{\bigotimes}\limits_{t\in T}\al H._t$ with
the inner product ${(\iota(\b x.),\,\iota(\b y.))}:=\hbox{``$\prod\limits_{t\in T}
(x_t,\, y_t)_t$''}$ whenever the right hand side makes sense.
There are many possibilities, but here we recall the tensor product constructions of von Neumann~\cite{VN38}.
Let
\[
\al L.:=\Big\{\b x.\in\prod_{t\in T}\al H._t\,\Big|\sum_{t\in T}\big|\|x_t\|_t-1\big|<\infty\Big\}
\]
where we interpret the convergence of a sum (resp.\ product)
over an uncountable set $T$ as convergence of the net of finite
partial sums, resp., products.
For sums such as $S:=\sum\limits_{t\in T} \alpha_t,$ $\alpha_t\in\C,$
this implies that
only countably many summands ${\{\alpha_{t_n}\mid n\in\N\}}$ are non-zero and that $S=
\sum\limits_{n=1}^\infty\alpha_{t_n},$ and it  converges absolutely (cf. Lemmas~2.3.2 and 2.3.3
in~\cite{VN38}).
Moreover, we have that $P=\prod\limits_{t\in T}|\alpha_t|<\infty$
if and only if either $\alpha_t=0$ for some $t$ (in which case $P=0$), or else
$\sum\limits_{t\in T}\big||\alpha_t|-1\big|<\infty$
(cf.~\cite[Lemma~2.4.1]{VN38}). We will not need to use general products
 $P=\prod\limits_{t\in T} \alpha_t,$ $\alpha_t\in\C,$ for which the convergence is a more difficult
 notion (cf. Lemma~2.4.2 and Definition~2.5.1
in~\cite{VN38}).

Thus $\b x.\in\al L.$ implies that $\|x_t\|_t=1$ for all $t\in T\backslash R$ where
$R$ is at most countable, and that the product $\prod\limits_{t\in T}\|x_t\|_t$
converges.
Obviously, any $\b x.$ such that $\|x_t\|_t=1$ for all $t\in T$ is in $\al L.\,.$
Note that if $\b x.\in\al L.$ then $\eqc{\b x.}.\subset\al L.$ also.
For $\bx, \by \in \cL$, we define
\begin{equation}
\label{VNeqc}
\b x.\approx\b y. \qquad\hbox{if}\qquad
\sum_{t\in T}\big|(x_t,y_t)_t-1\big|<\infty\,.
\end{equation}
Then $\approx$ is an equivalence relation by Lemma~3.3.3 in~\cite{VN38}, and we denote its
equivalence classes by $\eqcc{\b x.}.$. Observe that if $\b x.\in\al L.$ then
$\eqc{\b x.}.\subset\eqcc{\b x.}.\,,$ and moreover, each $\mathord{\approx}\hbox{-equivalence}$
class contains an $\b a.\in\al L.$ such that $\|a_t\|_t=1$ for all $t\in T$ (cf. Lemma~3.3.7
in~\cite{VN38}).

\begin{defn}\label{def:2.12}
Given such an $\b a.\in\eqcc{\b x.}.\subset\al L.\,,$
we can define an inner product on $\br{\b a.}.$ by sesqui-linear extension of
\[
\big(\iota(\b x.),\,\iota(\b y.)\big):=\prod\limits_{t\in T}
(x_t,\, y_t)_t\quad\hbox{for}\quad \b x.\sim\b a.\sim\b y.\,.
\]
(Note that the infinite products occurring here have only finitely many entries
different from $1$ hence are unproblematic).
Denote the closure of $\br{\b a.}.$ w.r.t. this Hilbert norm by
 $\mathop{\bigotimes}\limits_{t\in T}\!\!{}^{[\b a.]}\al H._t$. Then
this is von Neumann's ``incomplete direct product,'' and it contains
${\rm Span}\big(\iota(\eqcc{\b a.}.)\big)$
as a dense subspace (cf. Lemma~4.1.2 in~\cite{VN38}).
The direct sum of the spaces $\mathop{\bigotimes}\limits_{t\in T}\!\!{}^{[\b a.]}\al H._t$ where
we take one representative $\b a.$ from each  $\mathord{\approx}\hbox{-equivalence}$ class,
is von Neumann's ``complete direct product'' (cf. Lemma~4.1.1 in~\cite{VN38}).
An analogous associativity theorem to Theorem~\ref{Assoc} holds for this
complete direct product (cf. Theorem~VII in~\cite{VN38}).
\end{defn}

Next, consider the case where $(\cA_t)_{t \in T}$ is a
family of $*$-algebras, each equipped with a bounded Hilbert space $*$-representation
$\pi_t:\al A._t\to\al B.(\al H._t)\,.$
For any $\b A.\in\prod\limits_{t\in T}\al A._t$ we can define a linear map
$\pi\big(\iota\s{\al A.}.(\b A.)\big)$ on $\mathop{\otimes}\limits_{t\in T}\al H._t$ by
\[
\pi\big(\iota\s{\al A.}.(\b A.)\big)\iota(\b x.)=\mathop{\otimes}\limits_{t\in T}\pi_t(A_t)x_t
=\iota\big(\pi(\b A.)\b x.\big)\quad\hbox{for all}\quad\b x.\in\prod_{t\in T}\al H._t
\]
where 
$\big(\pi(\b A.)\b x.\big)_t:=\pi_t(A_t)x_t$ for all $t\in T$.
Then $\pi$ is a representation, because
it is one for each entry. To obtain Hilbert space $*$-representations from $\pi$, we need
to restrict it to suitable pre-Hilbert subspaces of $\mathop{\bigotimes}\limits_{t\in T}\al H._t$
hence need to restrict to those $\b A.$ such that $\pi\big(\iota\s{\al A.}.(\b A.)\big)$ preserves the Hilbert
space involved (and produces a bounded operator).

\begin{defn}
Consider the Hilbert space completion
$\mathop{\bigotimes}\limits_{t\in T}\!\!{}^{[\b a.]}\al H._t$
of $\br{\b a.}.\,,$ as above.
When the algebras $\al A._t$ are all unital, then $\br{\b 1.}.\subset
\mathop{\bigotimes}\limits_{t\in T}\al A._t$ is a $*$-subalgebra, where
$(\b 1.)_t=\1_t\in\al A._t$ for all $t\in T\,.$
Then $\pi(\b A.)\b x.\in\eqc{\b a.}.$ for all $\b x.\in\eqc{\b a.}.\subset\prod\limits_{t\in T}\al H._t$ and $\b A.\sim\b 1.\,.$
In particular, $\pi\big(\iota\s{\al A.}.(\b A.)\big)$ preserves $\br{\b a.}.$
and it is bounded, since it is a tensor product of a finite tensor product (of bounded operators)
with the identity operator. Thus it extends to a bounded operator on
 $\mathop{\bigotimes}\limits_{t\in T}\!\!{}^{[\b a.]}\al H._t\,.$
This defines a $*$-representation of the $*$-algebra $\br{\b 1.}.$ on the (stabilized) tensor product
$\mathop{\bigotimes}\limits_{t\in T}\!\!{}^{[\b a.]}\al H._t\,,$ and it is the most commonly used
definition of a tensor representation.

When the $*$-algebras $\al A._t$ are not unital,
consider the case where they contain nontrivial hermitian projections
$P_t\in\cA_t$. Then, for any choice of such projections
$\b P.\in\prod\limits_{t\in T}\al A._t$,
the subspace $\br{\b P.}.
\subset\mathop{\bigotimes}\limits_{t\in T}\al A._t$ is a $*$-subalgebra.
For any $\ba \in \prod\limits_{t \in T} \cH_t$ with
$\pi_t(P_t)a_t = a_t$ for all $t \in T$, we can now
define a tensor product representation of $\br{\b P.}.$
on $\mathop{\bigotimes}\limits_{t\in T}\!\!{}^{[\b a.]}\al H._t$.
Below we will consider more general tensor product representations.
\end{defn}

\section{Semi--host algebras for Gaussians}
\label{GausH}

In this section, $\mu$ will be a fixed Gaussian product measure on $\R^\N$
and $\mu_n$ denotes its projection on the $n^{\rm th}$ component.
For $\bx \in \R^\N$ and $\by \in \R^{(\N)}$, we write
$\la \bx, \by \ra :=\sum\limits_{i=1}^\infty x_iy_i$ for the standard
pairing. Recall that from  $\mu$
one constructs a unitary representation
\[
\pi_\mu:\R^{(\N)}\to \al U.\big(L^2(\R^\N,\mu)\big)\quad\hbox{by}\quad
\big(\pi_\mu(\b x.)f\big)(\b y.):=\exp\big(i\bbrk\b x.,{\b y.}.\big)\,f(\b y.), \quad \b x.\in\R^{(\N)}, \;\b y.\in\R^\N.
\]
Then there is a unitary map $U:\mathop{\bigotimes}\limits_{n=1}^\infty\!\!{}^{[\b e.]}\al H._n
\to L^2(\R^\N,\,\mu)$, where $\al H._n:=L^2(\R,\,\mu_n)$.
The sequence $\b e.=(e_1,\,e_2,\ldots)$ of stabilizing vectors $e_n\in\al H._n$
is given by the constant functions $e_n(x)=1$ for all $x\in\R\,.$
Explicitly, $U$ is given by
\[
U(f_1\otimes f_2\otimes\cdots\otimes f_k\otimes e_{k+1}\otimes e_{k+2}\otimes\cdots)
(x_1,x_2,\ldots)=f_1(x_1)\cdot f_2(x_2)\cdots f_k(x_k)
\]
which is clearly a cylinder function on $\R^\N\,.$
Then $\pi_\mu=U\Big(\bigotimes\limits_{n=1}^\infty\pi_{\mu_n}\Big)U^{-1}$, where each
\[
\pi_{\mu_n}:\R\to \al U.\big(L^2(\R,\,\mu_n)\big)\quad\hbox{is}\quad
\big(\pi_{\mu_n}(x)f\big)(y):=e^{ixy}\,f(y)\quad \mbox{ for } \quad x,y\in\R\,.
\]
The stabilizing sequence defines a cyclic vector $\Omega:=
\mathop{\otimes}\limits_{n=1}^\infty e_n\,.$
Immediate calculation establishes that the corresponding
positive definite function satisfies:
\begin{equation}
\label{MuDecomp}
 \omega_\mu(\b t.):=\big(\Omega,\,\pi_\mu(\b t.)\Omega\big)=\int_{\R^\N}\exp\big(i\bbrk\b t.,{\b y.}.\big)\,d\mu(\b y.)
\quad \mbox{ for } \quad \bt \in\R^{(\N)}\,,
\end{equation}
which is part of the Bochner--Minlos Theorem (cf.\ \cite{GV64}).
We will show that it expresses the decomposition of a state into the pure states of a
(semi-) host algebra for $\R^{(\N)}\,,$ and that there is a similar expression for
other states (which is also part of the Bochner--Minlos theorem).

Specialize the notation of the last
section by setting: $T=\N$ and $X_t=C_0(\R)\cong C^*(\R)$ for all $t\,.$
We first try to define an appropriate infinite tensor
product $C^*$-algebra of all the
$C_0(\R)\hbox{'s}\,,$ which seems to be a problem because $C_0(\R)$ is nonunital,
and has no nontrivial projection. By the last section we always have the algebraic
tensor product $\mathop{\bigotimes}\limits_{k=1}^\infty C_0(\R)$,
but this is too large.
We want to look at its $*$-subalgebras of the type defined in Theorem~\ref{ProdCls}(iii),
and will consider the following multiplicative semigroups in $C_0(\R)$. For each $n\in\N\,,$ define
\[
V_n:=\big\{f\in C_0(\R)\;\big|\;f(\R)\subseteq[0,1],\;\;f\rest[-n,n]=1,\;\;
{\rm supp}(f)\subseteq[-n-1,\,n+1]\big\}
\]
and observe that it is a semigroup, that $\|f\|=1$ for all $f\in V_n$ and that
any sequence\break
$\{u_n\in V_n\,\mid\,n\in\N\}$ is an approximate identity for $C_0(\R)\,.$
Moreover $V_n\cdot V_m=V_n$ if $m>n$ and hence $\bigcup\limits_{n=1}^\infty V_n$ is a semigroup.
For each $f\in V_n$ we have the subsemigroup
$$V_n(f):=\{f^k\,\mid\,k\in\N\}\subset V_n,$$
and for these we also have that
$V_n(f)\cdot V_m(g)=V_n(f)$ if $m>n\,.$

For any sequence $\b f.=(f_1,\,f_2,\,\ldots)\in C_0(\R)^\N$
with $f_n\in V_{k_n}$ for all $n\,,$  we consider the $*$-algebra generated in
$\mathop{\bigotimes}\limits_{k=1}^\infty C_0(\R)$ by $\br{\b f.}.\,,$ and note
that
\begin{equation}
  \label{eq:1}
 \hbox{$*$-alg}\big(\br{\b f.}.\big)
={\rm Span}\big\{\br{\b f.^k}.\,\mid\, k \in \N \}
\subset\mathop{\bigotimes}_{i=1}^\infty C_0(\R), \quad\hbox{where}\quad   \big(\b f.^k\big)_n
:=f_n^k\;\forall\,n
\end{equation}
and for the equality we needed the fact that $C_0(\R)\cdot C_0(\R)=C_0(\R)$
(Remark~\ref{rem:2.11}),
 and Theorem~\ref{ProdCls}(i).

Next, we want to define a convenient representation of $\hbox{$*$-alg}\big(\br{\b f.}.\big)$
to provide us with a $C^*$-norm to close it in.
We will show that there are $\b f.$ for which we can define a representation
of $\hbox{$*$-alg}\big(\br{\b f.}.\big)$ on $\mathop{\bigotimes}\limits_{n=1}^\infty\!\!{}^{[\b e.]}\al H._n$
in a natural way.

\begin{pro}
\label{RepLm}  We now have:
\begin{itemize}
\item[\rm(i)] Let $P_k$ denote multiplication of functions on $\R$ by
$\chi\s{[-k,k]}.\,.$ Then there exists  a sequence $(k_i)_{i\in\N}$ such that
${\sum\limits_{n=1}^\infty\big|(P_{k_n}e_n,\,e_n)_n-1\big|}<\infty\,.$
\item[\rm(ii)] Fix a sequence $(k_i)_{i\in\N}$ as in (i) as well as
$\b f. \in \prod\limits_{j = 1}^\infty V_{k_j}$. Then there
is a $*$-representation $\pi\s{\b e.}.:\hbox{$*$-alg}\big(\br{\b f.}.\big)
\to\al B.\Big(\mathop{\bigotimes}\limits_{n=1}^\infty\!\!{}^{[\b e.]}\al H._n\Big)$
such that
\[
\pi\s{\b e.}.\Big(\mathop{\otimes}_{n=1}^\infty g_n\Big)\mathop{\otimes}_{k=1}^\infty c_k
=\mathop{\otimes}_{n=1}^\infty g_n c_n\in
{\textstyle\mathop{\bigotimes}\limits_{n=1}^\infty\!\!{}^{[\b e.]}}\al H._n
\]
for all  $\b g.\sim\b f.^\ell,$ $\b c.\sim\b e.$ and $\ell\in\N$,
and where $g_n c_n$ is  the usual pointwise
product of functions on $\R\,.$
\end{itemize}
\end{pro}
\begin{beweis}
(i) For any $\varepsilon>0$, there is a $k\in\N$ such that
${\big|(P_ke_n,\,e_n)_n-1\big|}<\varepsilon$ by the Monotone Convergence Theorem.
Thus there is  a sequence $(k_i)_{i\in\N}$ such that
${\sum\limits_{n=1}^\infty\big|(P_{k_n}e_n,\,e_n)_n-1\big|}<\infty\,.$ \chop
(ii) Recall from Definition~\ref{def:2.12} that ${\rm Span}\big(\iota(\eqcc{\b e.}.)\big)$ is dense in the closure
$\mathop{\bigotimes}\limits_{n=1}^\infty\!\!{}^{[\b e.]}\al H._n$
of $\br{\b e.}.\,,$ where
\[
 \eqcc{\b e.}.=\Big\{\b v.\in\prod_{n=1}^\infty\al H._n\;\Big|\;
\sum_{n=1}^\infty\big|\|v_n\|_n-1\big|<\infty\quad\hbox{and}\quad
\sum_{n=1}^\infty\big|(e_n,\,v_n)_n-1\big|<\infty\Big\}\,.
\]
With the given choice of $(k_i)_{i\in\N}$ and $\b f.$ we have
\begin{eqnarray*}
 (P_{k_n}e_n,\,e_n)_n=\mu_n\big([-k_n,\,k_n]\big)&
\leq &\int_{-k_n-1}^{k_n+1}f_n(x)\,d\mu_n(x)
=(f_ne_n,\,e_n)_n \leq 1 \\[1mm]
\hbox{so that}\quad\qquad\big|(f_ne_n,\,e_n)_n-1\big|&\leq &\big|(P_{k_n}e_n,\,e_n)_n-1\big|,
\end{eqnarray*}
and hence $\sum\limits_{n=1}^\infty\big|(f_ne_n,\,e_n)_n-1\big|<\infty\,.$
As $(f_j)^\ell\in V_{k_j}$ for all $\ell\in\N\,,$ we have in fact that
 ${\sum\limits_{n=1}^\infty\big|(f_n^\ell e_n,\,e_n)_n-1\big|}<\infty$ for all $\ell\in\N\,.$
This implies that ${\sum\limits_{n=1}^\infty\big|\|f_n^\ell e_n\|_n^2-1\big|}<\infty$
which implies via Lemma~3.3.2 in~\cite{VN38} that
${\sum\limits_{n=1}^\infty\big|\|f_n^\ell e_n\|_n-1\big|}<\infty\,.$
Hence $\b f.^\ell\cdot\b e.\in\eqcc{\b e.}.$ and so
\[
 \Big(\mathop{\otimes}_{n=1}^\infty f_n^\ell\Big)\Big(\mathop{\otimes}_{k=1}^\infty e_k\Big)
=\mathop{\otimes}_{n=1}^\infty f_n^\ell e_n\in
{\textstyle\mathop{\bigotimes}\limits_{n=1}^\infty\!\!{}^{[\b e.]}}\al H._n\,.
\]
Since any $\b c.\sim\b e.$ only differs from $\b e.$ in finitely many entries, the
convergence arguments above will still hold if we replace $\b e.$ by $\b c.\,.$
Likewise, we can replace $\b f.^\ell$ by any $\b g.\sim\b f.^\ell\,,$ i.e.,
we have shown that
$\mathop{\otimes}\limits_{n=1}^\infty g_n c_n\in
{\mathop{\bigotimes}\limits_{n=1}^\infty\!\!{}^{[\b e.]}}\al H._n$
for all $\b g.\sim\b f.^\ell$ and $\b c.\sim\b e.\,.$
Since the multiplication map
$$ \Big(\bigcup\limits_{\ell\in\N}\eqc{\b f.^\ell}.\Big)\times\eqc{\b e.}.
\to \bigotimes_{n = 1}^\infty \cH_n, \quad
(\b g.,\b c.)\mapsto \mathop{\otimes}\limits_{n=1}^\infty g_n c_n $$
is multilinear, it defines a bilinear map on
${\rm Span}\Big(\bigcup\limits_{\ell\in\N}\br{\b f.^\ell}.\Big)\times\br{\b e.}.\,,$ denoted by $(a,b) \mapsto \pi_{\b e.}(a)b$,
thus obtaining the formula for $\pi_{\b e.}$ in
the theorem. That $\pi_{\b e.}$ is a representation of $\hbox{$*$-alg}\big(\br{\b f.}.\big)$ follows from the explicit formula,
and the $*$-property is also clear. It remains to show that each
$\pi_{\b e.}(a)$ is bounded (hence extends as a bounded operator to
${\mathop{\bigotimes}\limits_{n=1}^\infty\!\!{}^{[\b e.]}}\al H._n$).
It suffices to check this for the generating elements
of $\hbox{$*$-alg}\big(\br{\b f.}.\big)$.
Let $a\in\br{\b f.}.$ with $\ba \sim \b f.$:
\[
a= (a_1 \otimes \cdots \otimes a_p) \otimes f_{p+1}\otimes f_{p+2}\otimes\cdots
\]
for some $p<\infty\,.$ Moreover any $b\in\br{\b e.}.$ can
also be written in the form:
\[
b= b_p \otimes e_{p+1}\otimes e_{p+2}\otimes\cdots
\quad \mbox{ with } \quad
b_p \in \bigotimes_{j = 1}^p \cH_j, \]
where we may take the same $p$ as in the preceding expression (e.g. by adjusting the initial part).
Then
$$ \|\pi_{\b e.}(a)b\| = \|A_p b_p\|
\cdot\prod_{k=p+1}^\infty\|f_ke_k\|, \quad \mbox{ where } \quad
A_p v = (a_1 \otimes \cdots \otimes a_p)v. $$
Since $A_p$ is bounded on the completion $\mathop{\hat{\bigotimes}}\limits_{j=1,\ldots, p} \al H._j$  of
$\bigotimes\limits_{j=1}^p \al H._j,$ we have
$\|A_pb_p\|\leq\|A_p\|\cdot\|b_p\|\,,$ and as  $\|f_ke_k\|\leq\|e_k\|=1$,
we see that
\[
\|\pi_{\b e.}(a)b\|^2 \leq \|A_p\|^2\|b_p\|^2\cdot\prod_{k=p+1}^\infty\|e_k\|^2=
\|A_p\|^2\cdot\|b\|^2
\]
and hence $\pi_{\b e.}(a)$ is a bounded operator on $\br{\b e.}.$ so extends to
a bounded operator on ${\mathop{\bigotimes}\limits_{n=1}^\infty\!\!{}^{[\b e.]}}\al H._n\,.$
\end{beweis}

\begin{defn}
Thus for any $\bff \in \prod\limits_{j =1}^\infty V_{k_j}$, we can define
$$\al L._\mu[\b f.]:=C^*\left(\pi\s{\b e.}.\big(\hbox{$*$-alg}\big(\br{\b f.}.\big)\big)\right)
\subset\al B.\Big({\mathop{\bigotimes}\limits_{n=1}^\infty\!\!{}^{[\b e.]}}\al H._n\Big)\,.$$

\end{defn}

\begin{rem}
Recall that we also have the unitaries $\pi_\mu(\R^\N)\subset
\al U.\big({\mathop{\bigotimes}\limits_{n=1}^\infty
\!\!{}^{[\b e.]}}\al H._n\big),$
where
\begin{eqnarray*}
\pi_\mu(\b x.)\mathop{\otimes}\limits_{k=1}^\infty c_k
&=&\mathop{\otimes}\limits_{n=1}^\infty \big(\exp_{x_n}\!\!\cdot c_n\big)\in\br{\b e.}.,\;\;\;
\b x.\in\R^{(\N)}\,,\;\;\b c.\sim\b e.,\;\;\exp_{x_n}(t):=e^{ix_nt}\,.
\end{eqnarray*}
Then
$$ \pi_\mu(\b x.)\cdot\pi\s{\b e.}.\big(\iota(\b g.)\big)
=\pi\s{\b e.}.\big(\iota(\b g.)\big)\cdot\pi_\mu(\b x.)
=\pi\s{\b e.}.\Big(\mathop{\otimes}\limits_{n=1}^\infty \big(\exp_{x_n}\!\!\cdot g_n \big)\Big)\in\al L._\mu[\b f.], $$
for all  $\b x.\in\R^{(\N)},$ $\b g.\sim\b f.^\ell$  and $\ell\in\N\,.$
The inclusion needed the fact that $\b x.$ has only finitely many nonzero entries,
and that$\exp_{x_n}\!\!\cdot C_0(\R)\subset C_0(\R).$
Thus $\pi_\mu(\R^{(\N)})\cdot\al L._\mu[\b f.] \subset \cL_\mu[\b f.]$.
Since for each $\b x.\in\R^{(\N)}$
we can find a sequence $\big(\iota(\b g._n)\big)_{n\in\Z}\subset\br{\b f.}.$
such that $\pi_{\b e.}\big(\iota(\b g._n)\big)\cdot\pi_\mu(\b x.)$ converges in norm
to $\pi_\mu(\b x.),$ we have a faithful embedding of
$\R^{(\N)}$ as unitaries into the multiplier algebra
$M\big(\al L._\mu[\b f.]\big)$ denoted $\eta \: \R^{(\N)} \to M\big(\al L._\mu[\b f.]\big).$
In the next section we will investigate to what extent $\al L._\mu[\b f.]$ is  a host
algebra of $\R^{(\N)}\,.$
\end{rem}

\begin{lem}
\label{SepChar} With $\b f.$ as in {\rm Proposition~\ref{RepLm}(ii)},
we have
\begin{itemize}
\item[(i)] The $C^*$-algebra $\al L._\mu[\b f.]$ is separable.
\item[(ii)] Let $\omega$ be a pure state on $\al L._\mu[\b f.],$
and let $\wt{\omega}$ be its strict extension to
the unitaries \chop $\eta(\R^{(\N)})\subset M\big(\al L._\mu[\b f.]\big)\,.$
Then $\wt{\omega}\circ\eta$ is a character and there exists an
element $\ba \in \R^\N$ with
$\wt{\omega}(\eta(\b x.))=\exp\big(i\bbrk\b x.,{\b a.}.\big)$
for all $\b x.\in\R^{(\N)}$.
\end{itemize}
\end{lem}
\begin{beweis}
(i) Since $\pi\s{\b e.}.\big(\hbox{$*$-alg}\big(\br{\b f.}.\big)\big)$ is dense
in $\al L._\mu[\b f.]$, where
\[   
\hbox{$*$-alg}\big(\br{\b f.}.\big)={\rm Span}\big\{\br{\b f.^k}.\,\mid\, 
k\in\N\big\} \qquad           
\hbox{and}\qquad\br{\b f.^k}.=\bigcup_{m=1}^\infty\Big\{\Big({\textstyle\bigotimes\limits_{\ell=1}^m C_0(\R)}\Big)
\otimes f_{m+1}^k\otimes f_{m+2}^k\otimes\cdots\Big\},
\]  
(i) follows immediately from the separability of $C_0(\R)$.

\noindent (ii) As $\al L._\mu[\b f.]$ is commutative, any pure state $\omega$ of it is a point evaluation, hence
a $*$-homomorphism. Thus the strict extension $\wt{\omega}$ to $\eta(\R^{(\N)})\subset M\big(\al L._\mu[\b f.]\big)$
is also a $*$-homomorphism, hence $\wt{\omega}\circ\eta$ is a character. The restriction of $\wt{\omega}\circ\eta$ to the subgroup
$\R^n\subset\R^{(\N)}$ is still a character, and it is continuous (since it is determined
by the factor $\mathop{\otimes}\limits_{j=1}^nC_0(\R)$ in $\al L._\mu[\b f.]$
which is the group algebra of $\R^n$) hence of the form
$\wt{\omega}\circ\eta(\b x.)=\exp(i\b x.\cdot\b a.^{(n)})$ for some $\b a.^{(n)}\in\R^n\,.$
Since $\wt{\omega}\circ\eta$ is a character on all of $\R^{(\N)},$ the family
${\{\b a.^{(n)}\in\R^n\,\mid\,n\in\N\}}$ is a consistent family, i.e.,
if $n<m$ then $\b a.^{(n)}$ is the first $n$ entries of $\b a.^{(m)}.$
Thus there is an $\b a.\in\R^\N$ such that $\b a.^{(n)}$ is the first $n$ entries of
$\b a.$ for any $n\in\N\,.$ Then
$\wt{\omega}\circ\eta(\b x.)=\exp\big(i\bbrk\b x.,{\b a.}.\big)$
since for any $\b x.\in\R^n\subset\R^{(\N)}$ this restricts to the previous formula
for $\wt{\omega}\circ\eta\,.$
%
\end{beweis}

Since $\al L._\mu[\b f.]$ is separable and commutative,
it follows from Theorem~II.2.2 in~\cite{Dav}
that all its cyclic representations are multiplicity free,
and hence by Theorem~4.9.4 in~\cite{Ped},
for any state $\omega$ on $\al L._\mu[\b f.]$,
there is a regular Borel probability measure $\nu$ on the states $\wp(\al L._\mu[\b f.])$ concentrated on
the pure states $\wp_p(\al L._\mu[\b f.])$ such that
\begin{equation}
  \label{eq:4}
 \omega(A)=\int_{\wp_p(\al L._\mu[\b f.])}\varphi(A)\,d\nu(\varphi)\quad\forall\,A\in\al L._\mu[\b f.]\,.
\end{equation}
We will show that this decomposition produces similar decompositions to
the one in \eqref{MuDecomp} for other continuous positive definite functions than $\omega_\mu\,.$

Since $\al L._\mu[\b f.]$ is separable, it has a countable approximate
identity $\{E_n\}_{n\in\N}\subset\al L._\mu[\b f.]$ (cf. Remark~3.1.1\cite{Mur}).
For a state $\omega$ on $\al L._\mu[\b f.],$ let $\wt{\omega}$ be its strict extension to
the unitaries $\eta(\R^{(\N)})\subset M\big(\al L._\mu[\b f.]\big)\,,$ then  we have
for any countable approximate
identity $\{E_n\}_{n\in\N}\subset\al L._\mu[\b f.]$ that
\begin{eqnarray*}
\wt{\omega}\circ\eta(\b x.)&=&\lim_{n\to\infty}\omega(\eta(\b x.) E_n)=
\lim_{n\to\infty}\int_{\wp_p(\al L._\mu[\b f.])}\varphi(\eta(\b x.) E_n)\,d\nu(\varphi) \\[1mm]
&=& \int_{\wp_p(\al L._\mu[\b f.])}\lim_{n\to\infty}\varphi(\eta(\b x.)  E_n)\,d\nu(\varphi)
=\int_{\wp_p(\al L._\mu[\b f.])}\wt{\varphi}\circ\eta(\b x.)\,d\nu(\varphi)
\end{eqnarray*}
where we used the Lebesgue dominated convergence theorem in the second line, since
${\big|\varphi(\eta(\b x.)  E_n)\big|}\leq1$ and the
constant function $1$ is integrable.

By Lemma~\ref{SepChar}(ii) we can define a map
\[
\xi:\wp_p(\al L._\mu[\b f.])\to\R^\N \qquad\hbox{by}\qquad \wt{\varphi}\circ\eta(\b x.)=\exp\big(i\bbrk\b x.,\xi(\varphi).\big)
\;\;\quad \mbox{ for } \quad \;\b x.\in\R^{(\N)}\,,
\]
so using $\xi$ we define a probability measure $\wt{\nu}$ on $\R^\N$ by
$\wt{\nu}:=\xi_*\nu,$ and so:
\begin{equation}
\label{OmDecomp}
 \wt{\omega}\circ\eta(\b x.)=\int_{\R^\N}\exp\big(i\bbrk\b x.,{\b y.}.\big)\,d\wt{\nu}(\b y.)
\quad\mbox{ for } \quad \,\b x.\in\R^{(\N)}\,,
\end{equation}
which generalises the integral representation \eqref{MuDecomp} to
those positive definite functions
$\wt\omega$ which are strict extensions of states of $\al L._\mu[\b f.]$
(and this includes $\omega_\mu).$
We will obtain the full Bochner--Minlos theorem for $\R^{(\N)}$ in a
$C^*$-algebraic context, if we can show that every continuous
normalized positive definite function is of this type for some $\mu$ and some~$\b f.\,.$
This is what we will do in the next section.

\section{Semi-host algebras for $\R^{(\N)}$}
\label{ParHRN}

Inspired by the good properties which we found for $\al L._\mu[\b f.]$
above, we now examine more general versions of these algebras.
The semi-host algebras which we obtain will be the building blocks
for the algebra hosting the full representation theory of $\R^{(\N)},$ which will be constructed in the next section.

For the rest of this section we fix a sequence
$(k_n)_{n \in \N} \in \N^\N$ and
$\b f.\in\prod\limits_{n=1}^\infty V_{k_n}$ such that
$\br{\b f.}.\not=0\,.$
Then we have that
\begin{eqnarray}
\label{StAlgF}
  \hbox{$*$-alg}\big(\br{\b f.}.\big)&=&{\rm Span}\big\{\br{\b f.^k}.\,\mid\,
k \in \N\big\} =\ilim\al A._m[\b f.], \quad \mbox{ where} \\[1mm]
\al A._m[\b f.]&:=&{\rm Span}\big\{A_1\otimes\cdots\otimes A_m\otimes
f_{m+1}^k\otimes f_{m+2}^k\otimes\cdots\,\mid\,A_i\in C_0(\R)\;\;\forall\; i\in\N,\;k\in\N\big\}\nonumber
\end{eqnarray}
and  the inductive limit is w.r.t. set inclusion of the $*$-algebras  $\al A._m[\b f.]\subset\al A._\ell[\b f.]$
if $m<\ell\,.$ By the Associativity Theorem~\ref{Assoc},
we can write
\[
\al A._m[\b f.]=\Big({\textstyle\mathop{\bigotimes}\limits_{k=1}^m}C_0(\R)\Big)\otimes
\Big(\hbox{$*$-alg}\big({\textstyle\mathop{\bigotimes}\limits_{j=m+1}^\infty}f_j\big)\Big). \]
The natural $C^*$-norm on on the first factor is clear, but not on the second factor.
So we next investigate possible bounded $*$-representations to provide $\hbox{$*$-alg}\big(\br{\b f.}.\big)$
with a $C^*$--norm. Since $\hbox{$*$-alg}\Big(\mathop{\bigotimes}\limits_{j=m+1}^\infty f_j\Big)$
is generated by the single element
$E:=\mathop{\bigotimes}\limits_{j=m+1}^\infty f_j$,
any representation $\pi$ of this $*$-algebra is given by specifying the
single operator $\pi(E)\,.$ Since $E$ is positive, we require  $\pi(E)\geq 0\,,$ and as we
want a tensor norm on the larger $\hbox{$*$-alg}\big(\br{\b f.}.\big),$ we need that
$\|\pi(E)\|\leq\prod\limits_{j=m+1}^\infty \|f_j\|=1\,.$
\begin{lem}
\label{RepSet}
Let $\bff \in \prod\limits_{n \in \N} V_{k_n}$ and let
$\big\{\pi_k:C_0(\R)\to\al B.(\al H.)\,\mid\,k\in\N\big\}$ be a set
of $*$-representations on the same space with commuting ranges. Then
\begin{itemize}
\item[(i)] The strong limit $F_k^{(\ell)}:=\slim\limits_{n\to\infty}\pi_k(f_k^\ell)\cdots
\pi_n(f_n^\ell)\in\al B.(\al H.)$ exists, and
$0\leq F_k^{(\ell)}\leq \1$ for $k, \ell \in \N$.
\item[(ii)] $P[\b f.]:=\slim\limits_{k\to\infty}F_k^{(\ell)}$
(an increasing limit) is a projection independent of
$\ell\in\N$ satisfying $F_k^{(\ell)}P[\b f.]=F_k^{(\ell)}\,.$
\item[(iii)] Let $Q\in\al B.(\al H.)$ be such that
$0\leq Q\leq\1\,,$ and such that it commutes with
$\pi_k(C_0(\R))$ for each $k\in\N.$
Let $A:=A_1\otimes\cdots\otimes A_m\otimes f_{m+1}^\ell\otimes f_{m+2}^\ell\otimes\cdots
\in\hbox{$*$-alg}\big(\br{\b f.}.\big)$ and define
\[
 \pi_Q(A):=\pi_1(A_1)\,\pi_2(A_2)\cdots\pi_{m}(A_{m})\,F_{m+1}^{(\ell)}Q^{\ell}\,.
\]
Then $\pi_Q$ defines a $*$-representation $\pi_Q:\hbox{$*$-alg}\big(\br{\b f.}.\big)\to\al B.(\al H.)\,.$
\item[(iv)]
The representation  $\pi_Q$ is non-degenerate if and only if
all $\pi_i$ are non-degenerate, $P[\b f.]=\1$
and $\ker Q=\{0\}\,.$
If $\pi_Q$ is degenerate, $\ker Q = 0,$
  and all $\pi_j$ are non-degenerate,
then  $P[\b f.]$ is the projection onto the essential subspace of $\pi_Q.$
\end{itemize}
\end{lem}
\begin{beweis}
 (i) Since the operators $\pi_k(f_k^\ell),\;\pi_j(f_j^\ell)\in\al B.(\al H.)$ commute and are positive,
it follows from joint spectral theory that their product $\pi_k(f_k^\ell)\cdot\pi_j(f_j^\ell)$ is also a
positive operator. From  $\pi_k(f_k^\ell)\leq\1$ for
all $k,\,\ell\in\N,$ we derive that
$\pi_k(f_k^\ell)\cdot\pi_j(f_j^\ell)\leq\pi_k(f_k^\ell)$ and hence,
for a fixed $k$, the operators
$C_n:=\pi_k(f_k^\ell)\cdots\pi_n(f_n^\ell)$ form
a decreasing sequence of commuting positive operators.
Thus, by~\cite[Thm~4.1.1, p.~113]{Mur},
$C_n$ converges in the strong operator topology to some limit
$F_k^{(\ell)}\,.$ It is clear that $F_k^{(\ell)}$ is positive, and using
\[
 \|T\|=\sup\big\{\big|(\psi,\,T\psi)\big|\;\mid\,\psi\in\al H.,\;\|\psi\|=1\big\}
\qquad\hbox{whenever}\qquad T=T^*\,,
\]
it follows from $\|C_n\|=\big\|\pi_k(f_k^\ell)\cdots\pi_n(f_n^\ell)\big\|\leq 1$ for all $n$ that
$\|F_k^{(\ell)}\|\leq 1$ and hence that $0\leq F_k^{(\ell)}\leq \1\,.$

(ii) By definition, $F_k^{(\ell)}=\pi_k(f_k^\ell)F_{k+1}^{(\ell)}$ and $0\leq\pi_k(f_k^\ell)\leq\1$ and so
the commuting sequence of operators $\big(F_k^{(\ell)}\big)_{k\in \N}$
is increasing, and bounded above
by $\1\,.$ Thus it follows again from Theorem~4.1.1 in~\cite{Mur} that the
strong limit $P^{(\ell)}[\b f.]:=\slim\limits_{k\to\infty}F_k^{(\ell)}$ exists, is positive and
bounded above by $\1\,.$ Since the operator product is jointly strong operator continuous
on bounded sets, we get
\begin{eqnarray*}
 F_k^{(\ell)}P^{(\ell)}[\b f.]
&=&\slim_{n\to\infty}\pi_k(f_k^\ell)\cdots\pi_{n-1}(f_{n-1}^\ell)
\cdot\slim_{n\to\infty} F_n^{(\ell)}\\[1mm]
&=&\slim_{n\to\infty}\pi_k(f_k^\ell)\cdots\pi_{n-1}(f_{n-1}^\ell)\, F_n^{(\ell)} =
\slim_{n\to\infty}F_k^{(\ell)}=F_k^{(\ell)}\,.
\end{eqnarray*}
Thus by $P^{(\ell)}[\b f.]=\slim\limits_{k\to\infty}F_k^{(\ell)}=\slim\limits_{k\to\infty}F_k^{(\ell)}
P^{(\ell)}[\b f.]=\left(P^{(\ell)}[\b f.]\right)^2$ and the fact that $P^{(\ell)}[\b f.]$ is positive
we conclude that it is a projection. To see that $ P^{(\ell)}[\b f.]$ is independent of $\ell\,,$
note that for $k\leq m$ we have:
\begin{eqnarray}
 F_k^{(\ell)}F_m^{(j)}&=&\slim_{n\to\infty}\pi_k(f_k^\ell)\cdots\pi_n(f_n^\ell)
\cdot\slim_{p\to\infty} \pi_m(f_m^j)\cdots\pi_p(f_p^j)\nonumber\\[1mm]
&=&\slim_{n\to\infty}\pi_k(f_k^\ell)\cdots\pi_{m-1}(f_{m-1}^\ell)\,\pi_m(f_m^{\ell+j})
\cdots\pi_n(f_n^{\ell+j})\nonumber\\[1mm]
\label{piFj}
&=&\pi_k(f_k^\ell)\cdots\pi_{m-1}(f_{m-1}^\ell)\,F_m^{(\ell+j)}.
\end{eqnarray}
This leads to
$$P^{(\ell)}[\b f.]\cdot P^{(j)}[\b f.]= \slim_{k\to\infty}F_k^{(\ell)}\slim_{m\to\infty}F_m^{(j)}
 =\slim_{n\to\infty} F_n^{(\ell)} F_n^{(j)}
= \slim_{n\to\infty}F_n^{(\ell+j)}=P^{(\ell+j)}[\b f.].$$
However, each $P^{(\ell)}[\b f.]$ is idempotent, i.e.,
$P^{(\ell)}[\b f.]=P^{(2\ell)}[\b f.]$
for all $\ell\in\N\,,$ hence $P^{(\ell)}[\b f.]$ is independent of $\ell\,.$

(iii) Since $\hbox{$*$-alg}\big(\br{\b f.}.\big)=\ilim\al A._m[\b f.]=\bigcup\limits_{m\in\N}\al A._m[\b f.]\,,$
it suffices to show that $\pi_Q$ defines a $*$-representation on each $*$-algebra $\al A._m[\b f.]\,,$
and that $\pi_Q$ restricts to its correct values on any  $\al A._k[\b f.]\subset\al A._m[\b f.]$ for $k<m\,.$
Recall that
\[
\al A._m[\b f.]=\Big({\textstyle\mathop{\bigotimes}\limits_{k=0}^m}C_0(\R)\Big)\otimes
\Big(\hbox{$*$-alg}\big({\textstyle\mathop{\bigotimes}\limits_{j=m+1}^\infty}f_j\big)\Big)\,.
\]
Now
$$\pi_a^{(m)}:\mathop{\bigotimes}\limits_{k=0}^mC_0(\R)\to\al B.(\al H.), \quad
\pi_a^{(m)}(A_1\otimes\cdots\otimes A_m):=\pi_1(A_1)\cdots\pi_m(A_m)$$
is a well-defined $*$-representation
obtained by the universal property of the tensor product.
Moreover, since $\hbox{$*$-alg}\Big(\mathop{\bigotimes}\limits_{j=m+1}^\infty f_j\Big)$ is generated by
a single element not satisfying any polynomial relation,
the assignment
$\pi_b^{(m)}\Big(\mathop{\bigotimes}\limits_{j=m+1}^\infty f_j\Big):=F_{m+1}^{(1)}Q\geq 0$
defines a $*$-representation $\pi_b^{(m)}:\hbox{$*$-alg}\Big(\mathop{\bigotimes}\limits_{j=m+1}^\infty f_j\Big)
\to\al B.(\al H.)\,.$ Note from Equation~(\ref{piFj}) that $F_{m+1}^{(k)}\cdot F_{m+1}^{(\ell)}=F_{m+1}^{(k+\ell)}\,,$
which leads to the factorization
\[
\pi_Q\big(A_1\otimes\cdots\otimes A_m\otimes f_{m+1}^\ell\otimes f_{m+2}^\ell\otimes\cdots\big)
=\pi_a^{(m)}(A_1\otimes\cdots\otimes A_m)\cdot
\pi_b^{(m)}\Big(\big({\textstyle\mathop{\bigotimes}\limits_{j=m+1}^\infty} f_j\big)^\ell\Big)\,.
\]
Thus, since it is multilinear, we obtain a
linear map $\pi_Q$
on $\al A._m[\b f.]\,,$ and as the ranges of the $*$-representations $\pi_a$ and $\pi_b$ commute, $\pi_Q$ is a $*$-representation on $\al A._m[\b f.]\,.$
For $k<m$ we have from the definition that
$$ \pi_b^{(k)}\Big(\mathop{\bigotimes}\limits_{j=k+1}^\infty f_j\Big)=F_{k+1}^{(1)}Q=
\pi_{k+1}(f_{k+1})\cdots\pi_m(f_m)F_{m+1}^{(1)}Q $$
and hence
\begin{eqnarray*}
&&\pi_a^{(m)}(A_1\otimes\cdots\otimes A_k\otimes f_{k+1}\otimes\cdots f_m)\cdot
\pi_b^{(m)}\Big({\textstyle\mathop{\bigotimes}\limits_{j=m+1}^\infty} f_j\Big)
=\pi_a^{(k)}(A_1\otimes\cdots\otimes A_k)\cdot
\pi_b^{(k)}\Big({\textstyle\mathop{\bigotimes}\limits_{j=k+1}^\infty} f_j\Big)
\end{eqnarray*}
so it is clear that the value of $\pi_Q$ on $\al A._k[\b f.]\subset\al A._m[\b f.]$
is the same as the restriction of the map $\pi_Q$ defined on $\al A._m[\b f.].$
Hence $\pi_Q$ is consistently defined as a $*$-representation of
$\hbox{$*$-alg}\big(\br{\b f.}.\big)\,.$

(iv) Note that by $F_k^{(\ell)}P[\b f.]=F_k^{(\ell)}\,,$ we have  $\pi_Q(A)P[\b f.]=\pi_Q(A)$
for all $A\in\hbox{$*$-alg}\big(\br{\b f.}.\big)\,,$ hence,
if $P[\b f.]\not=\1$, then
 $\pi_Q\left(\hbox{*-alg}\big(\br{\b f.}.\big)\right)$ has null spaces, i.e., $\pi_Q$ is degenerate.
 Likewise, if $\ker Q\not=\{0\}$ then $\pi_Q$ is degenerate.
 Moreover, if any $\pi_i$ is degenerate, then since by commutativity:
 \[
 \pi_Q\big(A_1\otimes\cdots\otimes A_m\otimes f_{m+1}^\ell\otimes f_{m+2}^\ell\otimes\cdots\big)
 =\pi_1(A_1)\cdots\widehat{\pi_i(A_i)}\cdots\pi_{m}(A_{m})\,F_{m+1}^{(\ell)}Q^\ell\pi_i(A_i)\,,
 \]
 where the hat means omission, it follows that $\pi_Q$ is also degenerate.

 Conversely, let $\pi_Q$ be degenerate, i.e.,
there is a nonzero $\psi\in\al H.$ such that
 $\pi_Q(A)\psi=0$ for all $A\,,$ hence
 \[
 \pi_Q\big(A_1\otimes\cdots\otimes A_m\otimes f_{m+1}^\ell\otimes f_{m+2}^\ell\otimes\cdots\big)\psi
 =\pi_1(A_1)\cdots\pi_{m}(A_{m})\,F_{m+1}^{(\ell)}Q^\ell\psi=0
 \]
for all $A_i\in C_0(\R)$ and $m,\,\ell\in\N\,.$
 If all $\pi_j$ are non-degenerate, then it follows inductively
that  $F_{m}^{(\ell)}Q^\ell\psi=0$ for all $m$ and $\ell$.
If $\ker Q = 0$, then $F_{m}^{(\ell)}\psi=0$ for all $m$, hence
  $P[\b f.]\psi=0$, i.e., $P[\b f.]\not=\1\,.$

  By the last step we also see that when $\pi_Q$ is degenerate, $\ker Q = 0,$
  and all $\pi_j$ are non-degenerate, then
  $P[\b f.]$ is zero on the null space of $\pi_Q.$
  Since $F_k^{(\ell)}P[\b f.]=F_k^{(\ell)}$  by (ii) it follows from the definition
  of $\pi_Q$ that $\pi_Q(A)P[\b f.]=\pi_Q(A)$ for all $A\in\hbox{*-alg}\big(\br{\b f.}.\big)\,.$
  Thus $P[\b f.]$ is the identity on the essential subspace of $\pi_Q,$ i.e. it
  is the projection onto this essential subspace.
\end{beweis}

\begin{defn}
Using this lemma, we can now investigate natural representations of $\hbox{$*$-alg}\big(\br{\b f.}.\big)\,.$
Start with the universal representation of $\R^{(\N)}$ denoted
$\pi_u:\R^{(\N)}\to \al U.(\al H._u)$ which we recall, is the direct sum of the cyclic strong--operator
continuous unitary  representations of $\R^{(\N)},$ one from each unitary equivalence class.
Since for the $k^{\rm th}$component we have an inclusion
$\R\subset\R^{(\N)}$ by $x\longrightarrow{(0,\ldots,0,x,0,0,\ldots)}$ ($k^{\rm th}$entry), $\pi_u$
restricts to a representation on the $k^{\rm th}$component, denoted by
$\pi_u^k:\R\to \al U.(\al H._u)\,.$ By the host algebra property of $C^*(\R)\cong C_0(\R),$ this
produces a unique representation $\pi_u^k:C_0(\R)\to \al B.(\al H._u)\,,$ which is non-degenerate.  Since the set of representations
${\big\{\pi_u^k:C_0(\R)\to \al B.(\al H._u)\,\mid\,k\in\N\big\}}$ have commuting ranges,
we can apply Lemma~\ref{RepSet}, with $Q=\1$, to define a representation
$\pi_u:\hbox{$*$-alg}\big(\br{\b f.}.\big)\to\al B.(\al H._u)$ by an abuse of notation.
Below we will use the notation $F_{u,k}^{(\ell)}$ for the operator $F_{k}^{(\ell)}$ of $\pi_u.$
\end{defn}

\begin{defi}
\label{DefLf}
The $C^*$-algebra $\al L.[\b f.]$ is the $C^*$-completion of $\pi_u\big(\hbox{$*$-alg}\big(\br{\b f.}.\big)\big)$ in $\al B.(\al H._u)\,.$
\end{defi}
\begin{rem}
\begin{itemize}
\item[(1)]
We see directly from (\ref{StAlgF}) and the separability of $C_0(\R)$ that $\al L.[\b f.]$ is separable.
\item[(2)]
Observe that the representation
$\pi_u$ of $*\mbox{-alg}(\brf)$
may be degenerate. Although all $\pi^u_k$  are non-degenerate, it is possible that
$P[\b f.]\not=\1\,.$  By Lemma~\ref{RepSet}(iv) it then follows that $P[\b f.]$ is the projection
onto the essential subspace of $\pi_u.$
\item[(3)]
Since $\al L.[\b f.]\subset\al B.(\al H._u)$ is given as a concrete C*-algebra, this selects the class
of those representations of $\al L.[\b f.]$ which are normal maps w.r.t. the $\sigma\hbox{--strong}$ topology of
$\al B.(\al H._u)$ on $\al L.[\b f.].$ We will say that such a representation $\pi$ is normal w.r.t. the defining representation
$\pi_u.$ This will be the case if the vector states of $\pi\left( \al L.[\b f.]\right)$ are normal states for
$\pi_u\left( \al L.[\b f.]\right)$ (cf. Proposition~7.1.15 \cite{KR2}).
\item[(4)]
From Fell's Theorem~\cite[Thm.~1.2]{Fell} we know that
any state of
$\al L.[\b f.]$ is in the weak-*-closure of the convex hull of
the vector states of $\pi_u$.
\end{itemize}
\end{rem}
We will need the following proposition.
\begin{pro}
\label{Lfactor}
If $S\subset\N$ is  a finite subset, then
\begin{itemize}
\item[(i)] there is a  $C^*$-algebra $\al B._S[\b f.]\subset\al B.(\al H._u)$ and a copy of
the C*-complete tensor product
$\al L.^S:=\mathop{\hat\bigotimes}\limits_{s\in S} C_0(\R)$ in $\al B.(\al H._u)$
such that
\[
\al L.[\b f.]=C^*\big(\al L.^S\cdot\al B._S[\b f.]\big)
\cong \al L.^S\hat\otimes\al B._S[\b f.]\,.
\]
\item[(ii)]
the natural embeddings $\zeta_S:M(\al L.^S)\to M\big(\al L.[\b f.]\big)
=M\big(\al L.^S\hat\otimes\al B._S[\b f.]\big)$ by $\zeta_S(M)(A\otimes B):=
(M\cdot A)\otimes B$ for all $A\in\al L.^S$ and $B\in\al B._S[\b f.]$
are topological embeddings w.r.t. the strict topology on each bounded subset
of $M(\al L.^S)\,.$ Moreover, $\al L.^S$ is dense in $M(\al L.^S)$ w.r.t. the relative
strict topology of $M\big(\al L.[\b f.]\big)\,.$
\item[(iii)] The group homomorphism
$\eta \: \R^{(\N)} \to M\big(\al L.[\b f.]\big)$ is strictly continuous.
\end{itemize}
\end{pro}
\begin{beweis}
(i) By associativity (Theorem~\ref{Assoc}): $\bigotimes\limits_{k=1}^\infty C_0(\R)
=\Big(\bigotimes\limits_{s\in S} C_0(\R)\Big)\otimes\Big(\bigotimes\limits_{t\in \N\backslash S} C_0(\R)\Big)\,,$
and so, applying this to $\hbox{$*$-alg}\big(\br{\b f.}.\big),$ and using the fact that it is the span of
elementary tensors of the type
${A_1\otimes\cdots\otimes A_m\otimes f^\ell_{m+1}\otimes f^\ell_{m+2}\otimes\cdots}$ with $A_i\in C_0(\R)$ and $m,\,\ell\in\N,$
we get
\[
 \hbox{$*$-alg}\big(\br{\b f.}.\big)=
\Big({\textstyle\bigotimes\limits_{s\in S}} C_0(\R)\Big)\otimes\Big(\hbox{$*$-alg}\big(\br{\b f._{\N\backslash S}}.\big)
\Big)\,,
\]
where  $(\b f._{\N\backslash S})_t = f_t$ for $t\in\N\backslash S$
and
$\hbox{$*$-alg}\big(\br{\b f._{\N\backslash S}}.\big)$ denotes the $*$-algebra generated in
$\bigotimes\limits_{t\in \N\backslash S} C_0(\R)$ by
\[
 \Big\{\mathop{\otimes}_{t\in \N\backslash S}g_t\,\Big|\,\b g.\in\mathop{\prod}_{t\in \N\backslash S}C_0(\R),\;\;
\b g.\sim\b f._{\N\backslash S}\Big\}\,.
\]
Below, we need unital algebras, so adjoin identities, and define
\[
 \al C._0:=\Big(\C\1+{\textstyle\bigotimes\limits_{s\in S}} C_0(\R)\Big)\otimes\Big(\C\1+
\hbox{$*$-alg}\big(\br{\b f._{\N\backslash S}}.\big)\Big)\subset{\textstyle\bigotimes\limits_{k=1}^\infty}\big(
\C\1+ C_0(\R)\big)
\]
which contains $\hbox{$*$-alg}\big(\br{\b f.}.\big)$ as a $*$-ideal. Since $\pi_u\big(\hbox{$*$-alg}\big(\br{\b f.}.\big)\big)$
acts non-degenerately on its essential space $\al H._{\rm ess}\subset\al H._u,$ it determines a unique extension of $\pi_u$ to a representation \break
$\pi_u:\al C._0\to\al B.(\al H._u)\,,$ if we let the null space of $\pi_u$ be
$\al H._{\rm ess}^\perp.$ Define $\al C.:={C^*\big(\pi_u(\al C._0)\big)}
={C^*\big(\al A.\cdot\al B.\big)}$ where
\[
 \al A.:=C^*\Big(\pi_u\Big(\big(\C\1+{\textstyle\bigotimes\limits_{s\in S}} C_0(\R)\big)\otimes\1\Big)\Big)
\qquad\hbox{and}\qquad
\al B.:=C^*\Big(\pi_u\Big(\1\otimes\big(\C\1+
\hbox{$*$-alg}\big(\br{\b f._{\N\backslash S}}.\big)\Big)\Big)\,.
\]
Thus the unital $C^*$-algebra $\al C.$ is generated by the
two commmuting unital $C^*$-algebras
$\al A.$ and~$\al B..$ Moreover, since $\pi_u$ contains tensor representations (w.r.t. the two factors
of $\hbox{$*$-alg}\big(\br{\b f.}.\big)$ above), it follows that if $AB=0$ for an $A\in\al A.$ and a
$B\in\al B.,$ then either $A=0$ or $B=0\,.$ Thus
by \cite[Ex.~2, p.~220]{Tak}, it follows that
 $\al C.\cong\al A.\hat\otimes\al B.$, where the tensor $C^*$-norm is unique, since both  $\al A.$ and $\al B.$
are commutative, hence nuclear. We conclude that
the original $C^*$-norm defined on $\al C.$ is in fact a
cross--norm. Since its restriction to
\[
 \hbox{$*$-alg}\big(\br{\b f.}.\big)=
\Big({\textstyle\bigotimes\limits_{s\in S}} C_0(\R)\Big)\otimes\Big(\hbox{$*$-alg}\big(\br{\b f._{\N\backslash S}}.\big)
\Big)\subset\al C._0
\]
is still a cross--norm, and the latter is unique by commutativity of the algebras (given the norms on the factors),
 it follows from
$C^*\big[\pi_u\big(\bigotimes\limits_{s\in S} C_0(\R)\big)\big]
=\hat\bigotimes_{s\in S}C_0(\R)$    that
\begin{eqnarray*}
 \al L.[\b f.]&=&\Big({\textstyle\bigotimes\limits_{s\in S}}
C_0(\R)\Big)\otimes
C^*\Big[\pi_u\Big(\1\otimes\big(\hbox{$*$-alg}\big(\br{\b f._{\N\backslash S}}.\big)\Big)\Big] =
\al L.^S\hat\otimes\al B._S[\b f.]\\[1mm]
&=&C^*\Big[\pi_u\Big(\big({\textstyle\bigotimes\limits_{s\in S}} C_0(\R)\big)\otimes\1\Big)
\cdot\pi_u\Big(\1\otimes\big(\hbox{$*$-alg}\big(\br{\b f._{\N\backslash S}}.\big)\Big)\Big],
\end{eqnarray*}
where $\al B._S[\b f.]:=C^*\Big[\pi_u\Big(\1\otimes
\big(\hbox{$*$-alg}\big(\br{\b f._{\N\backslash S}}.\big)\big)\Big)\Big]\,.$
\chop
(ii)  This follows from (i) and Lemma~A.2 in~\cite{GrNe}.

\noindent (iii) Since $\eta(\R^{(\N)})$ consists of unitary multipliers,
it suffices to verify that the set of all elements $A \in \cL[\bff]$
for which the map
$$ \eta^A\: \R^{(\N)} \to \cL[\bff], \quad \bx \to \eta(\bx)A$$
is continuous span a dense subalgebra. To establish this, let
$A = \iota(\by)$ for some $\by \sim \bff^k$ for some $k \in \N$.
Now $\R^{(\N)}$ is a topological direct limit, so that
it suffices to verify continuity on the finite dimensional subgroups
$\R^n.$ For these, it follows from the strict continuity of the action  of the group
$\R^n$ on its $C^*$-algebra $C^*(\R^n) \cong C_0(\R^n)$ and the fact that
by part~(i) we have
$$ \cL[\bff] \cong C_0(\R^n) \hat\otimes \cA, $$
for a $C^*$-algebra $\cA$, where $\R^n$ acts by unitary multipliers on the first tensor factor
and the identity on the second factor.
\end{beweis}

Note that for $S=\{1,2,\ldots,n\}\,,$ the map $\zeta_S$ identifies $\R^n\subset UM(\al L.^S)$ with
the unitaries $\R^n\subset\R^{(\N)}\subset  UM(\al L.[\b f.])\,.$
Below we will abbreviate the notation to $\al L.^{(n)}:=\al L.^{\{1,2,\ldots,n\}}=\hat\bigotimes_{k=1}^nC_0(\R)\,.$
For ease of notation, we sometimes also omit explicit indication of the embeddings $\zeta_S,$ using inclusions instead.

Next, let $\pi:\al L.[\b f.]\to\al B.(\al H._\pi)$ be a given fixed non-degenerate $*$-representation.
Let $\wt\pi$ denote the strict extension of $\pi$ to $M(\al L.[\b f.])\,,$
so that
$\pi_k:=\wt{\pi}\restriction\al L.^{\{k\}}$ and $\pi^{(n)}:=\wt{\pi}\restriction\al L.^{(n)}$ are the strict
extensions of $\pi$ to $\al L.^{\{k\}}\subset M(\al L.^{\{k\}})\xhookrightarrow{\zeta_{\{k\}}}
M(\al L.[\b f.])$ and $\al L.^{(n)}\subset M(\al L.^{(n)})\xhookrightarrow{\zeta_{\{1,\ldots,n\}}}
M(\al L.[\b f.])$ respectively.
Then ${\big\{\pi_k\,\mid\,k\in\N\big\}}$ is a set of non-degenerate representations
with commuting ranges as in Lemma~\ref{RepSet}, hence we specialize its notation to:
\[
F_{\pi,k}^{(\ell)}:=\slim\limits_{n\to\infty}\pi_k(f_k^\ell)\cdots
\pi_n(f_n^\ell)\in\al B.(\al H._\pi)\qquad\hbox{and}\qquad
 P_\pi[\b f.]:=\slim\limits_{k\to\infty}F_{\pi,k}^{(\ell)}\in\al B.(\al H._\pi)\,.
\]
Since the commuting sequence of operators $\big(F_{\pi,k}^{(\ell)}\big)_{k=1}^\infty$ is increasing,
$P_\pi[\b f.]\not=\1$ implies that there is a nonzero $\psi\in\al H._\pi$
such that $F_{\pi,k}^{(\ell)}\psi=0$ for all $k$ and $\ell\,.$

We will show in the next proposition that, for a certain choice of $Q$,
there is a  representation $\pi_Q$ constructed as in Lemma~\ref{RepSet}
from the set ${\big\{\pi_k\,\mid\,k\in\N\big\}}$ which coincides with $\pi\,.$
\begin{pro}
\label{repLf}
Fix a non-degenerate $*$-representation $\pi:\al L.[\b f.]\to\al B.(\al H._\pi)$
with $\al H._\pi\not=\{0\}.$
\begin{itemize}
\item[(i)] Let $B_n:=\wt{\pi}
\big(\mathop{\overbrace{\1\otimes\cdots\otimes\1}}\limits^{n-1\;{\rm factors}}
\otimes f_n\otimes f_{n+1}\otimes\cdots\big)$.
Then the strong limit
$Q := \slim\limits_{n\to\infty}B_n$ exists and satisfies $0< Q\leq\1.$
\item[(ii)] If $A:=A_1\otimes\cdots\otimes A_m\otimes f_{m+1}^\ell\otimes f_{m+2}^\ell\otimes\cdots
\in\hbox{$*$-alg}\big(\br{\b f.}.\big),$ then
\[
 \pi(A)=\pi_1(A_1)\,\pi_2(A_2)\cdots\pi_{m}(A_{m})\,F_{\pi, m+1}^{(\ell)}\,Q^\ell=\pi_Q(A),
\]
i.e.,  $\pi_Q = \pi\restriction{*\mbox{-alg}(\brf)}$.
Moreover  $P_\pi[\b f.]=\1$ and $\ker Q=\{0\}\,.$
\item[(iii)]
Let $\pi^{(n)}:\al L.^{(n)}\to\al B.(\al H._\pi)$ denote
the strict extension of $\pi$ to
$\al L.^{(n)} \subeq M(\cL(\brf))$. Then
\[
\pi(L_1\otimes L_2\otimes\cdots)=\slim_{n\to\infty}\pi^{(n)}\big(L_1\otimes L_2\otimes\cdots\otimes L_n\big)Q^\ell
\]
for all elementary tensors $L_1\otimes L_2\otimes\cdots\in\br{\b f.^\ell}.\subset\hbox{$*$-alg}\big(\br{\b f.}.\big).$
\end{itemize}
\end{pro}
\begin{beweis}
(i) We need to prove this claim in greater
 generality than stated above, for use in the subsequent part.
By definition, we have for $A:=A_1\otimes\cdots\otimes A_m\otimes f_{m+1}^\ell\otimes f_{m+2}^\ell\otimes\cdots
\in\hbox{$*$-alg}\big(\br{\b f.}.\big),$ that
\[
 \pi_u(A)=\pi_u^1(A_1)\,\pi_u^2(A_2)\cdots\pi_u^{m}(A_{m})\,F_{u,m+1}^{(\ell)}\in\al L.[\b f.],
\]
where $F_{u,k}^{(\ell)}:=\slim\limits_{n\to\infty}\pi_u^k(f_k^\ell)\cdots
\pi_u^n(f_n^\ell)=\wt{\pi}_u\big(\1\otimes\cdots\otimes\1\otimes f^\ell_n\otimes
f^\ell_{n+1}\otimes\cdots\big)\in\al B.(\al H._u)$. Hence we have that
$F_{u,n}^{(\ell)}\in M\big(\al L.[\b f.]\big)$.
Thus the operator
\[
B_n^{(\ell)}:=\wt{\pi}\big(\mathop{\overbrace{\1\otimes\cdots\otimes\1}}^{n-1\;{\rm factors}}\otimes f^\ell_n\otimes f^\ell_{n+1}\otimes\cdots\big)=\wt{\pi}\big(F_{u,n}^{(\ell)}\big)
\]
satisfies $0\leq B_n^{(\ell)}\leq\1$ since $0\leq F_{u,n}^{(\ell)}\leq \1\,.$
As $B_n^{(\ell)}=\pi_n(f_n^\ell)B_{n+1}^{(\ell)}$ and
$\pi_n(f_n^\ell)\leq\1$ is a positive operator commuting with $B_{n+1}^{(\ell)},$ we see that
$B_n^{(\ell)}\leq B_{n+1}^{(\ell)}.$ Thus the
strong limit $Q^{(\ell)}:= \slim\limits_{n\to\infty}B_n^{(\ell)}$ exists
by Theorem~4.1.1 in~\cite{Mur}, and satisfies $0< Q^{(\ell)}\leq\1$
 (note that $Q^{(\ell)}\not=0$ since $\pi$ is
non-degenerate and $\cH_\pi \not= \{0\}).$
Since the operator product is jointly strongly
continuous on bounded sets we have:
\begin{eqnarray*}
Q^{(\ell)}Q^{(m)}&=&\slim_{n\to\infty} \wt{\pi}\big(\1\otimes\cdots\otimes\1\otimes f^\ell_n\otimes f^\ell_{n+1}\otimes\cdots\big)
\slim_{k\to\infty} \wt{\pi}\big(\1\otimes\cdots\otimes\1\otimes f^m_k\otimes f^m_{k+1}\otimes\cdots\big)
\\[1mm]
&=&\slim_{n\to\infty} \wt{\pi}\big(\1\otimes\cdots\otimes\1\otimes f^\ell_n\otimes f^\ell_{n+1}\otimes\cdots\big)
\, \wt{\pi}\big(\1\otimes\cdots\otimes\1\otimes f^m_n\otimes f^m_{n+1}\otimes\cdots\big)
\\[1mm]
&=&\slim_{n\to\infty} \wt{\pi}\big(\1\otimes\cdots\otimes\1\otimes f^{\ell+m}_n\otimes f^{\ell+m}_{n+1}\otimes\cdots\big)= Q^{(\ell+m)}\,.
\end{eqnarray*}
Thus $Q^{(\ell)}=Q^\ell$ where $Q:=Q^{(1)}\,.$

(ii) Now $B_n^{(\ell)}=\wt{\pi}\big(\mathop{\overbrace{\1\otimes\cdots\otimes\1}}\limits^{n-1\;
{\rm factors}}\otimes f^\ell_n\otimes f^\ell_{n+1}\otimes\cdots\big)$
\begin{eqnarray}
&=& \wt{\pi}\big(\1\otimes\cdots\otimes\1\otimes f^\ell_n\otimes \1\otimes\cdots\big)\cdot
\wt{\pi}\big(\1\otimes\cdots\otimes\1\otimes f^\ell_{n+1}\otimes f^\ell_{n+2}\otimes\cdots\big)\nonumber\\[1mm]
&=&\pi_n(f_n^\ell)\,\wt{\pi}\big(\1\otimes\cdots\otimes\1\otimes f^\ell_{n+1}\otimes f^\ell_{n+2}\otimes\cdots\big)\nonumber\\[1mm]
&=&\slim_{k\to\infty} \pi_n(f_n^\ell)\cdots\pi_k(f_k^\ell)\,
\wt{\pi}\big(\mathop{\overbrace{\1\otimes\cdots\otimes\1}}^{k\;{\rm factors}}\otimes f^\ell_{k+1}\otimes f^\ell_{k+2}\otimes\cdots\big)\nonumber\\[1mm]
&=&\slim_{k\to\infty} \pi_n(f_n^\ell)\cdots\pi_k(f_k^\ell)\,\slim_{m\to\infty}
\wt{\pi}\big(\1\otimes\cdots\otimes\1\otimes f^\ell_{m+1}\otimes f^\ell_{m+2}\otimes\cdots\big)\nonumber\\[1mm]
\label{InterM1}
&=& F_{\pi,n}^{(\ell)}Q^{(\ell)}= F_{\pi,n}^{(\ell)}Q^\ell
\end{eqnarray}
where we used again the joint strong operator continuity of the product on bounded sets.
Let $A:=A_1\otimes\cdots\otimes A_m\otimes f_{m+1}^\ell\otimes f_{m+2}^\ell\otimes\cdots
\in\hbox{$*$-alg}\big(\br{\b f.}.\big).$ Then
\begin{eqnarray}
\pi(A)&=& \pi_1(A_1)\cdot\wt{\pi}\big(\1\otimes A_2\otimes A_3\otimes\cdots\otimes A_m\otimes f_{m+1}^\ell\otimes f_{m+2}^\ell\otimes\cdots\big)
=\cdots \nonumber\\[1mm]
&=&\pi_1(A_1)\,\pi_2(A_2)\cdots\pi_{m}(A_{m})\cdot
\wt{\pi}\big(\1\otimes \cdots\otimes \1\otimes f_{m+1}^\ell\otimes f_{m+2}^\ell\otimes\cdots\big) \nonumber\\[1mm]
\label{InterM2}
&=&\pi_1(A_1)\,\pi_2(A_2)\cdots\pi_{m}(A_{m})\cdot F_{\pi,m+1}^{(\ell)}Q^\ell=\pi_Q(A)
\end{eqnarray}
making use of (\ref{InterM1}) above.
Since $\pi$ is non-degenerate, it follows from Lemma~\ref{RepSet}(iii)
that $P_\pi[\b f.]=\1$ and $\ker Q=\{0\}\,.$
\chop
(iii) Note first that from Proposition~\ref{Lfactor}(ii) above
and \cite[Lemma~4.1 on p.203]{Tak} that
$\pi_1(A_1)\,\pi_2(A_2)\cdots\pi_{n}(A_{n})=\pi^{(n)}(A_1\otimes\cdots\otimes A_n)$ for all $A_i\in C_0(\R)\,.$
Thus, if we continue equation (\ref{InterM2}) above
\begin{eqnarray*}
\pi(A)&=& \pi_1(A_1)\,\pi_2(A_2)\cdots\pi_{m}(A_{m})\cdot F_{\pi,m+1}^{(\ell)}Q^\ell  \\[1mm]
&=&\pi_1(A_1)\,\pi_2(A_2)\cdots\pi_{m}(A_{m})\,\slim_{n\to\infty}\pi_{m+1}(f_{m+1}^\ell)\cdots
\pi_n(f_n^\ell)\,Q^\ell\\[1mm]
&=&\slim_{n\to\infty}\pi^{(n)}(A_1\otimes\cdots\otimes A_m\otimes f_{m+1}^\ell\otimes\cdots\otimes f_n^\ell)\, Q^\ell
\end{eqnarray*}
which establishes the claim.
\end{beweis}

\begin{defn}
Given a representation $\pi$ of $\al L.[\b f.],$ we will call its associated operator $Q$
its {\it excess}.
\end{defn}

This proposition creates a difficulty for the host algebra project, because by part~(iii) we can see that
to construct its representations, we need more information than what is contained in the representations
of $\R^{(\N)},$ i.e., we need the excess operators $Q\,.$ It is therefore very important to establish whether
there are representations $\pi_Q$ with $Q\not=\1$ (below we will see such $\pi_Q$ will
not be normal w.r.t. $\pi_u$).
\begin{pro}
\label{Qchar}
Let $\b f.$ be as before and let
$\big\{\pi_k:C_0(\R)\to\al B.(\al H.)\,\mid\,k\in\N\big\}$ be a set
of $*$-representations on the same space with commuting ranges. Then for any positive operator
$Q\in\al B.(\al H.)$ with $Q\leq\1$ which commutes with the ranges of all $\pi_k,$ we have that
$\pi_Q:\hbox{$*$-alg}\big(\br{\b f.}.\big)\to\al B.(\al H.)$ extends to a
$*$-representation of $\al L.[\b f.]\,.$
\end{pro}
\begin{beweis}
We show first that $\sigma\left(F_{u,k}\right)=[0,1].$
Let $\omega$ be a character of $\R^{(\N)}.$ Then since it is a one--dimensional subrepresentation
of $\pi_u$ there is a vector $\psi_\omega\in\al H._u$ such that ${\big(\psi_\omega,\pi_u(\b x.)\psi_\omega\big)}
=\omega(\b x.)$ for all $\b x.\in\R^{(\N)}.$
Then $\omega_k(h)= {\big(\psi_\omega,\pi^k_u(h)\psi_\omega\big)}$ for all $h\in\al L.^{\{k\}}=C_0(\R)$
is also a character, hence a point evaluation at a point $x_k^\omega\in\R$,
and in fact we obtain all point evaluations of $\al L.^{\{k\}}=C_0(\R)$ this way.
Thus
\[
F_{\omega,k}:=\slim_{n\to\infty}\omega_k(f_k)\cdots
\omega_n(f_n)=\lim_{n\to\infty}f_k(x_k^\omega)\cdots
f_n(x_n^\omega)=\prod_{n=k}^\infty f_n(x_n^\omega)\in[0,1]\,,
\]
and as we can choose our $\omega,$ hence points $x_k^\omega\in\R$ arbitrarily,
it is clear that we can find $\omega$ to set $F_{\omega,k}$ equal to any
value in $[0,1].$ Since
\[
F_{\omega,k}:=\lim_{n\to\infty}\omega_k(f_k)\cdots
\omega_n(f_n)=\wt{\omega}\big(\mathop{\overbrace{\1\otimes\cdots\otimes\1}}\limits^{k-1\;{\rm factors}}\otimes f^\ell_k\otimes f^\ell_{k+1}\otimes\cdots\big)=\big(\psi_\omega,F_{u,k}\psi_\omega)
\]
defines a character on ${\rm C}^*(F_{u,k})$ we see that $\sigma\left(F_{u,k}\right)=[0,1]\,.$
Since for
$\{\pi_k\,\mid\,k\in\N\}$ and $Q$ as in the initial hypotheses we always have
that $0\leq F_{\pi,k}Q\leq\1,$ it follows that
$\sigma\big(F_{\pi,k}Q\big)\subseteq[0,1]=\sigma\left(F_{u,k}\right)$ for all $k.$

Next, note that in a diagonalization of $F_{u,k}\geq 0$ we can write it
as $F_{u,k}(x)=x$ for $x\in\sigma(F_{u,k}),$ and hence
$\|p(F_{u,k})\|=\sup\big\{\big|p(x)\big|\,\mid\,x\in \sigma(F_{u,k})\big\}\,.$
From this it is immediate that $\sigma\left(F_{\pi,k}Q\right)\subseteq
\sigma\left(F_{u,k}\right)$ implies $\left\|p(F_{\pi,k}Q)\right\|\leq
\left\|p\big(F_{u,k}\big)\right\|$ for all  polynomials $p\,.$

Finally, recall that $\hbox{$*$-alg}\big(\br{\b f.}.\big)=\ilim\al A._m[\b f.]$ where
\[
\al A._m[\b f.]:={\rm Span}\big\{A_1\otimes\cdots\otimes A_m\otimes
f_{m+1}^k\otimes f_{m+2}^k\otimes\cdots\,\mid\,A_i\in C_0(\R)\;\;\forall\; i\in\N,\;k\in\N\big\}
\]
and the inductive limit is w.r.t. to the
inclusion $\al A._m[\b f.]\subset\al A._\ell[\b f.]\,.$
Thus $\al L.[\b f.]$ is the inductive limit of the
$C^*$-closures $\cL_m$ of $\pi_u\big(\al A._m[\b f.]\big)$
w.r.t. set inclusion. Since
\[
\al A._m[\b f.]=\Big({\textstyle\mathop{\bigotimes}\limits_{k=0}^m}C_0(\R)\Big)\otimes
\Big(\hbox{$*$-alg}\big({\textstyle\mathop{\bigotimes}\limits_{j=m+1}^\infty}f_j\big)\Big)\,,
\]
and the norm of $\al L.[\b f.]$ is a product norm by Proposition~\ref{Lfactor}(i), we have that
$\al L._m\cong\al L.^{(m)}\hat\otimes C^*(F_{u,m+1})\,.$
Next we define (as in
the proof of Lemma~\ref{RepSet}(iii)) two $*$-representations
$\pi_a^{(m)}:\mathop{\bigotimes}\limits_{k=0}^mC_0(\R)\to\al B.(\al H.)$
and $\pi_b^{(m)}:\hbox{$*$-alg}\Big(\mathop{\bigotimes}\limits_{j=m+1}^\infty f_j\Big)
\to\al B.(\al H.)$ as follows. First, we have that
$$\pi_a^{(m)}:\mathop{\bigotimes}\limits_{k=0}^mC_0(\R)\to\al B.(\al H.), \quad
\pi_a^{(m)}(A_1\otimes\cdots\otimes A_m):=\pi_1(A_1)\cdots\pi_m(A_m)$$
defines a well-defined $*$-representation by the universal property of the tensor product.
Moreover, since $\hbox{$*$-alg}\Big(\mathop{\bigotimes}\limits_{j=m+1}^\infty f_j\Big)$ is generated by
a single element not satisfying any polynomial relation,
the assignment
$\pi_b^{(m)}\Big(\mathop{\bigotimes}\limits_{j=m+1}^\infty f_j\Big):=F_{m+1}^{(1)}Q\geq 0$
defines a $*$-representation $\pi_b^{(m)}:\hbox{$*$-alg}\Big(\mathop{\bigotimes}\limits_{j=m+1}^\infty f_j\Big)
\to\al B.(\al H.)\,.$ Note from Equation~(\ref{piFj}) that $F_{m+1}^{(k)}\cdot F_{m+1}^{(\ell)}=F_{m+1}^{(k+\ell)}\,,$
which leads to the factorization
\[
\pi_Q\big(A_1\otimes\cdots\otimes A_m\otimes f_{m+1}^\ell\otimes f_{m+2}^\ell\otimes\cdots\big)
=\pi_a^{(m)}(A_1\otimes\cdots\otimes A_m)\cdot
\pi_b^{(m)}\Big(\big({\textstyle\mathop{\bigotimes}\limits_{j=m+1}^\infty} f_j\big)^\ell\Big)\,.
\]
Now $\pi_a^{(m)}$ has a unique extension to $\al L.^{(m)}\,,$ and as $\pi_b^{(m)}$ is defined
on the dense $*$-algebra $\hbox{$*$-alg}\Big(\mathop{\bigotimes}\limits_{j=m+1}^\infty f_j\Big)
={\big\{p\big(F_{u,k}\big)\,\mid\,p\;\;\hbox{a polynomial}\big\}}$ on which it is
continuous by the fact proven above, that $\big\|\pi_b^{(m)}\big(p(F_{u,k})\big)\big\|=\left\|p(F_{\pi,k}Q)\right\|\leq
\left\|p\big(F_{u,k}\big)\right\|.$ Thus
it extends uniquely to $C^*(F_{u,m+1}),$ hence
$\pi\s Q.$ has a unique continuous extension to $\al L._m\,.$ Since $\pi\s Q.$
respects the inductive limit structure (since it does so on the dense subalgebra
$\hbox{$*$-alg}\big(\br{\b f.}.\big)$ and is continuous on all $\al A._m$) it follows that
$\pi_Q$ extends uniquely to a continuous $*$-representation of $\al L.[\b f.]\,.$
\end{beweis}
We conclude that there is an abundance of representations $\pi$ of $\al L.[\b f.]$ with $Q\not=\1.$

Having investigated the representations of $\al L.[\b f.],$ we next consider its host algebra properties.
First label the unitary embedding $\eta:\R^{(\N)}\to M\big(\al L.[\b f.]\big)$ where
\begin{eqnarray*}
&&\!\!\!\!\!\!\!
 \eta(x_1,\ldots,x_n,0,0,\ldots)(L_1\otimes L_2\otimes\cdots)=\eta_1(x_1)L_1\otimes\cdots\otimes
\eta_n(x_n)L_n\otimes L_{n+1}\otimes L_{n+2}\otimes\cdots \\[1mm]
&&=\zeta\s{\{1,\ldots,n\}}.(x_1,\ldots,x_n)\big(L_1\otimes  L_2
\otimes \cdots\big)
\end{eqnarray*}
for all $(x_1,\ldots,x_n)\in\R^n,$ $n\in\N,$ $L_i\in\al L.^{\{i\}}=C_0(\R),$
and where $\eta_i:\R\to M({\rm C}^*(\R))$ is the usual unitary embedding.
Then the map $\eta^*:{\rm Rep}\big(\al L.[\b f.],\al H.\big)\to
{\rm Rep}\big(\R^{(\N)},\al H.\big)$ consists of the strict extension of
(non-degenerate) representations of $\al L.[\b f.]$ to $\eta(\R^{(\N)}),$ i.e.
\[
 \eta^*(\pi)(\b x.):=\slim_{\alpha\to\infty}\pi\big(\eta(\b x.)
E_\alpha\big)\quad \mbox{ for } \quad \bx \in \R^{(\N)} \]
and any approximate identity $\{E_\alpha\}\s\alpha\in\Lambda.$ in $\cL[\bff].$
Since $\al L.[\b f.]$ and $\R^{(\N)}$ are commutative, their irreducible representations are all
one-dimensional, hence $\eta^*$ takes irreducible representations to irreducible representations.
\begin{teo}
 \label{HostLf}
Given the preceding notation, we have that
\begin{itemize}
 \item[(i)] $\eta:\R^{(\N)}\to M\big(\al L.[\b f.]\big)$ is continuous w.r.t. the strict topology of
$M\big(\al L.[\b f.]\big)\,.$
\item[(ii)] Let ${\rm Rep}_0\big(\al L.[\b f.],\al H.\big)$ denote those non-degenerate $*$-representations
of $\al L.[\b f.]$ with excess operators $Q=\1$ (cf. Proposition~\ref{repLf}).
Then $\eta^*$ is injective on ${\rm Rep}_0\big(\al L.[\b f.],\al H.\big)\,.$
\item[(iii)] The range $\eta^*\left({\rm Rep}\big(\al L.[\b f.],\al H.\big)\right)$
is the same as $\eta^*\left({\rm Rep}_0\big(\al L.[\b f.],\al H.\big)\right)$
and consists of those $\pi\in{\rm Rep}\big(\R^{(\N)},\al H.\big)$ such that
$\1=\slim\limits_{k\to\infty}\wt{F}_k$ where
$\wt{F}_k:=\slim\limits_{n\to\infty}\pi_k(f_k)\cdots\pi_n(f_n)$ with
$\pi_k$ the unique representation in ${\rm Rep}\big(\al L.^{\{k\}},\al H.\big)$ such that $\eta_k^*(\pi_k)=\pi\restriction\R e_k$,
where $e_k \in \R^{(\N)}$ is the $k^{\rm th}$ basis vector.
\item[(iv)] For a state $\omega\in\wp\big(\al L.[\b f.]\big)\,,$ its GNS--representation $\pi^\omega$ is in
${\rm Rep}_0\big(\al L.[\b f.],\al H._\omega\big)$ if and only if
\[
 \omega\in\wp_0\big(\al L.[\b f.]\big):=\big\{\varphi\in\wp\big(\al L.[\b f.]\big)\,\mid\,
\lim_{n\to\infty}\wt\varphi\big(
\mathop{\overbrace{\1\otimes\cdots\otimes\1}}\limits^{n-1\;{\rm factors}}
\otimes f_n\otimes f_{n+1}\otimes\cdots\big)=1\big\}\,.
\]
Moreover, the restriction $\eta^*:\wp_0\big(\al L.[\b f.]\big)\to\wp\big(\R^{(\N)}\big)\equiv$states of $\R^{(\N)},$ is injective, with range consisting of
\[
\omega\in\wp\big(\R^{(\N)}\big)\quad\hbox{such that}\quad \lim_{k\to\infty}\lim_{n\to\infty}
\left(\Omega_\omega,\,\pi_k^\omega(f_k)\cdots\pi_n^\omega(f_n)\Omega_\omega\right)=1
\]
with $\pi_j^\omega$ as in (iii), and $\Omega_\omega$ is the cyclic GNS--vector.
\item[(v)]
$\pi$ is normal w.r.t. the defining representation
$\pi_u$ of $\al L.[\b f.]$ if and only if $Q=\1\,.$
\end{itemize}
\end{teo}
\begin{beweis}
(i) This is proven already in Proposition~\ref{Lfactor}(iii).

(ii) Let $\pi\in{\rm Rep}_0\big(\al L.[\b f.],\al H.\big)$ and let $\wt\pi$ be its strict extension
to $M\big(\al L.[\b f.]\big).$ As $\wt\pi$ is strictly continuous, (i) implies that the
unitary representation $\eta^*(\pi)=\wt{\pi}\circ\eta:\R^{(\N)}\to\al U.(\al H.)$ is strong operator
continuous. We need to show that $\eta^*$ is injective on ${\rm Rep}_0\big(\al L.[\b f.],\al H.\big).$
If  $\eta^*(\pi)=\eta^*(\pi')$ for two representations $\pi,\,\pi'\in{\rm Rep}_0\big(\al L.[\b f.],\al H.\big),$
then $\eta_{\{1,\ldots,n\}}^*(\pi)=\eta_{\{1,\ldots,n\}}^*(\pi')$ on $\R^n\subset\R^{(\N)}$
for all $n\in\N.$ But ${\rm Span}\big(\eta_{(n)}(\R^n)\big)\subset M(\al L.^{(n)})$ is strictly dense, and by
Proposition~\ref{Lfactor}(ii) this is still true for the strict topology of $M(\al L.[\b f.])\supset
\zeta\s{\{1,\ldots,n\}}.\big(M(\al L.^{(n)})\big).$ Thus $\wt{\pi}\restriction \zeta\s{\{1,\ldots,n\}}.\big(
\al L.^{(n)}\big)=\pi^{(n)}=\wt{\pi'}\restriction \zeta\s{\{1,\ldots,n\}}.\big(\al L.^{(n)}\big),$ i.e.,
$\pi$ and $\pi'$ produce the same
representation $\pi^{(n)}:\al L.^{(n)}\to\al B.(\al H.).$ Thus by Proposition~\ref{repLf}(iii)
(using $Q=\1$) we find
\[
\pi(L_1\otimes L_2\otimes\cdots)=\slim_{n\to\infty}\pi^{(n)}\big(L_1\otimes L_2\otimes\cdots\otimes L_n\big)
=\pi'(L_1\otimes L_2\otimes\cdots)
\]
for the elementary tensors in $\hbox{$*$-alg}\big(\br{\b f.}.\big)\,,$ i.e., $\pi=\pi'.$
Thus $\eta^*$ is injective on ${\rm Rep}_0\big(\al L.[\b f.],\al H.\big).$

(iii) To see that $\eta^*\left({\rm Rep}\big(\al L.[\b f.],\al H.\big)\right)
=\eta^*\left({\rm Rep}_0\big(\al L.[\b f.],\al H.\big)\right),$ note that for $\pi_Q$ as in Lemma~\ref{RepSet}:
\begin{eqnarray*}
&&\!\!\!\!\!\!\!
 \pi_Q\big(\eta(x_1,\ldots,x_n,0,0,\ldots)(L_1\otimes L_2\otimes\cdots)\big)\\
&&=\pi_Q\big(\eta_1(x_1)L_1\otimes\cdots\otimes
\eta_n(x_n)L_n\otimes L_{n+1}\otimes L_{n+2}\otimes\cdots \big)\\[1mm]
&&=\slim_{k\to\infty}\pi^{(k)}\big(\eta_1(x_1)L_1\otimes\cdots\otimes
\eta_n(x_n)L_n\otimes L_{n+1}\otimes L_{n+2}\otimes\cdots \otimes L_k\big)Q^\ell\qquad\hbox{by Prop.~\ref{repLf}(iii)}\\[1mm]
&&=\slim_{k\to\infty}\pi_1\big(\eta_1(x_1)L_1\big)\cdots\pi_n\big(
\eta_n(x_n)L_n\big)\pi_{n+1}\big(L_{n+1}\big)\pi_{n+2}\big(L_{n+2}\big)\cdots \pi_{k}\big(L_k\big)Q^\ell\\[1mm]
&&=\eta_1^*\pi_1(x_1)\cdots \eta_n^*\pi_n(x_n)\slim_{k\to\infty}\pi_1\big(L_1\big)\cdots\pi_{k}\big(L_k\big)Q^\ell\\[1mm]
&&=\eta_1^*\pi_1(x_1)\cdots \eta_n^*\pi_n(x_n)\,\pi_Q(L_1\otimes L_2\otimes\cdots)
\end{eqnarray*}
for $L_1\otimes L_2\otimes\cdots\in\br{\b f.^\ell}.\subset\hbox{$*$-alg}\big(\br{\b f.}.\big),$ which shows that
$$\eta^*(\pi_Q)(x_1,\ldots,x_n,0,0,\ldots)=\eta^*(\pi_\1)(x_1,\ldots,x_n,0,0,\ldots),$$ and establishes the claim.

To characterize the range of $\eta^*,$ let $\pi\in{\rm Rep}_0\big(\al L.[\b f.],\al H.\big)$ and note that as it is
non-degenerate, we have from Lemma~\ref{RepSet} that
\[
\1= P_\pi[\b f.]:=\slim\limits_{k\to\infty}F_{\pi,k}^{(\ell)}
\qquad\hbox{where}\qquad
F_{\pi,k}^{(\ell)}:=\slim\limits_{n\to\infty}\pi_k(f_k^\ell)\cdots
\pi_n(f_n^\ell)\in\al B.(\al H._\pi)\,,
\]
and $\pi_k=\wt{\pi}\restriction\al L.^{\{k\}}.$ From the uniqueness of the strict extension $\wt\pi$ on
$M(\al L.[\b f.])$ and the fact that the strict topology of $M(\al L.^{\{k\}})\subset M(\al L.[\b f.])$
coincides with that of $M(\al L.[\b f.])$ on bounded subsets,
we see that $\eta_k^*(\pi_k)=\eta^*\pi\restriction\R e_k$
and hence $\wt{F}_k=F_{\pi,k}^{(1)}.$ Thus $\1=\slim\limits_{k\to\infty}\wt{F}_k\,.$

Conversely, let $\pi\in{\rm Rep}\big(\R^{(\N)},\al H._\pi\big)$ be such that
$\1=\slim\limits_{k\to\infty}\wt{F}_k\,.$ We want to define
$\pi_{\al L.}\in{\rm Rep}_0\big(\al L.[\b f.],\al H._\pi\big)$
such that $\eta^*(\pi_{\al L.})=\pi\,.$ Consider first the case that $\pi$
is cyclic. Recall that $\al L.[\b f.]$ is the norm closure of
${\pi_u\big(\hbox{$*$-alg}\big(\br{\b f.}.\big)\big)}\,.$ By definition of $\pi_u,$ $\al H._\pi$ is a direct summand of
$\al H._u$
and there is a  projection $P_\pi\in{\pi_u(\R^{(\N)})'}$ such that
$\pi(\b x.)=P_\pi\pi_u(\b x.)\restriction\al H._\pi$ for all $\b x.\in\R^{(\N)}.$
 Then $\pi_k(A)=P_\pi\pi^k_u(A)\restriction\al H._\pi$ for all $A\in\al L.^{\{k\}},$
 and hence $\wt{F}_k= P_\pi F_{u,k}^{(1)}\restriction\al H._\pi\,.$
 We define $\pi_{\al L.}:\al L.[\b f.]\to\al B.(\al H._\pi)$ by $\pi_{\al L.}(A):=P_\pi\pi_u(A)\restriction\al H._\pi$
 which is obviously a $*$-representation, satisfying $F_{\pi_{\al L.},k}=\wt{F}_k,$ with excess $\1$
 (as it is normal w.r.t. $\pi_u$), and as
 \[
 P_{\pi_{\al L.}}[\b f.]=\slim_{k\to\infty}F_{\pi_{\al L.},k}=\slim_{k\to\infty}\wt{F}_k=\1
 \]
 by hypothesis, $\pi_{\al L.}$ is non-degenerate. Next, relax the requirement that $\pi$ be cyclic.
 Then $\pi$ is a direct sum of cyclic representations. Let
  $(\pi^c,\,\al H._c)$ be a cyclic subrepresentation of $\pi,$ and denote
 the projection onto $\al H._c$ by $P_c.$
 Since $\pi\restriction\R e_k$ also preserves $\al H._c$, it follows that
 $\pi_k^c(A)=P_c\pi^k_u(A)\restriction\al H._c$ for all $A\in\al L.^{\{k\}}.$
 Now, recalling that
 $\1=\slim\limits_{k\to\infty}\wt{F}_k$ where
$\wt{F}_k:=\slim\limits_{n\to\infty}\pi_k(f_k)\cdots\pi_n(f_n),$ we have that
\begin{eqnarray*}
\1\s{\al H._c}.=P_c\restriction\al H._c&=&P_c\,\slim_{k\to\infty}\wt{F}_k\restriction\al H._c
=\slim_{k\to\infty}\,\slim_{n\to\infty}\,P_c\pi_k(f_k)\cdots\pi_n(f_n)\restriction\al H._c  \\[1mm]
&=&\slim_{k\to\infty}\,\slim_{n\to\infty\,}\pi_k^c(f_k)\cdots\pi^c_n(f_n)=
\slim_{k\to\infty}\wt{F}_k^c
\end{eqnarray*}
where $\wt{F}_k^c:=\slim\limits_{n\to\infty}\pi_k^c(f_k)\cdots\pi_n^c(f_n).$
Thus, by the previous part we can construct a nondegenerate representation
$\pi_{\al L.}^c:\al L.[\b f.]\to\al B.(\al H._c)$ by $\pi_{\al L.}^c(A):=P_{\pi^c}\pi_u(A)\restriction\al H._{\pi_c}$
which is normal w.r.t. $\pi_u.$ Define $\pi_{\al L.}:\al L.[\b f.]\to\al B.(\al H._\pi)$
as the direct sum of all the $\pi_{\al L.}^c.$
Since this is normal w.r.t. $\pi_u$ and nondegenerate, we have that $\pi_{\al L.}\in{\rm Rep}_0\big(\al L.[\b f.],\al H._\pi\big).$

 Since the strict extension of $\pi_{\al L.}$ produces the same
 representations $\pi_k$ on $\al L.^{\{k\}}$ than obtained from $\pi\restriction\R e_k$, the strict
 extension of $\pi_{\al L.}$ must coincide on $\R^{(\N)}$ with $\pi,$ i.e.
 $\eta^*(\pi_{\al L.})=\pi\,.$

(iv) It is immediate from the definitions that if
$\pi^\omega\in{\rm Rep}_0\big(\al L.[\b f.],\al H._\omega\big)\,,$ then
$ \omega\in\wp_0\big(\al L.[\b f.]\big)\,.$ Conversely,
let $ \omega\in\wp_0\big(\al L.[\b f.]\big).$
Then, as $\al L.[\b f.]$ is commutative, we know $\al L.[\b f.]\cong C_0(X),$
with $X$ its spectrum. Then there is a probability measure $\mu$ on $X$ and a unitary
$U:\al H._\omega\to L^2(X,\mu)$ such that $(U\pi^\omega(h)\psi)(x) = h(x)\big(U\psi\big)(x)$ for all
$h\in C_0(X),$ $\psi\in\al H._\omega,$ $x\in X\,,$ and moreover $U\Omega_\omega=1\,.$
Then
\begin{eqnarray*}
 1&=& \lim_{n\to\infty}\wt\omega\big(
\mathop{\overbrace{\1\otimes\cdots\otimes\1}}\limits^{n-1\;{\rm factors}}
\otimes f_n\otimes f_{n+1}\otimes\cdots\big)= \left(\Omega_\omega,\,Q\Omega_\omega\right)\\[1mm]
&=& \int_X (UQU^{-1})(x)\,d\mu(x)\quad\hbox{and as}\quad 0<Q\leq\1
\quad\hbox{we have:} \\[1mm]
0 &=& \int_X \left|1-(UQU^{-1})(x)\right|\,d\mu(x).
\end{eqnarray*}
Hence $(UQU^{-1})(x)=1$ $\mu\hbox{--a.e.,}$ i.e., $Q=\1$ and thus
$\pi^\omega\in{\rm Rep}_0\big(\al L.[\b f.],\al H._\omega\big)\,.$

The last part of the claim now follows from this, (iii), and the observation that
$\eta^*\omega(g)=\left(\Omega_\omega,\,\eta^*\pi^\omega(g)\Omega_\omega\right)$
for all $g\in\R^{(\N)}.$ Note that the state condition on the range of $\eta^*$ implies
the operator condition in (iii) by a similar argument than the one above for $Q\,.$

(v) Let $\pi$ be  normal  w.r.t.~$\pi_u(\al L.[\b f.])$. Then it is continuous
on bounded sets w.r.t. the strong operator topologies
of both sides, hence
$Q=\slim\limits_{n\to\infty}B_n^{(1)}=\slim\limits_{n\to\infty}\wt{\pi}\big(F_{u,n}^{(1)}\big)
=\wt{\pi}\big(\slim\limits_{n\to\infty}F_{u,n}^{(1)}\big)\,.$
However, by Lemma~\ref{RepSet}(iv) we have that $P_u[\b f.]=\slim\limits_{n\to\infty}F_{u,n}^{(1)}$ is the projection
onto the essential subspace of $\pi_u(\al L.[\b f.]).$ Thus, since $\al L.[\b f.]$ is in fact defined in $\pi_u,$
it follows that $P_u[\b f.]$ is the identity for $\pi_u(\al L.[\b f.]),$
hence
$Q=\wt{\pi}\big(\slim\limits_{n\to\infty}F_{u,n}^{(1)}\big)=\1\,.$

Conversely, let $Q=\1,$ then by part (iii) $\eta^*\pi$ is a continuous representation
of $\R^{(\N)},$ and by
Proposition~\ref{repLf}(iii) (with $Q=\1)$ we have that
\[
\pi(L_1\otimes L_2\otimes\cdots)=\slim_{n\to\infty}\pi^{(n)}\big(L_1\otimes L_2\otimes\cdots\otimes L_n\big)
=\slim_{n\to\infty}\pi_1(L_1)\pi_2(L_2)\cdots\pi_n( L_n)
\]
for all elementary tensors $L_1\otimes L_2\otimes\cdots\in\br{\b f.^\ell}.\subset\hbox{$*$-alg}\big(\br{\b f.}.\big).$
This is precisely the formula in which Lemma~\ref{RepSet} defined representations on $\hbox{$*$-alg}\big(\br{\b f.}.\big)$
which we used to define $\pi_u.$ Now $\pi_u(\R^{(\N)})''={\{\pi_u^{(n)}(\al L.^{(n)})\,\mid\,n\in\N\}''}=
\pi_u(\al L.[\b f.])''$ and a similar equation holds for $\pi.$
Since the cyclic components of $\pi$ are contained in the direct summands of $\pi_u,$
there is a normal map $\varphi:\pi_u(\al L.[\b f.])''\to\al B.(\al H._\pi)$ such that
$\varphi\circ\pi_u=\pi.$ Thus $\pi$ is normal to $\pi_u.$
\end{beweis}

Thus, though $\al L.[\b f.]$ is not actually a host algebra for $\R^{(\N)},$ it does have
good properties, e.g., $\eta^*$ is bijective between two large sets of representations, and it takes
irreducible representations to irreducibles. In fact, using
the algebras $\al L.[\b f.]$,
we can now give a full $C^*$-algebraic interpretation of the
Bochner--Minlos Theorem.
Our aim is not to re--prove the Bochner--Minlos Theorem in the C*-context, but
just to identify the measures and decompositions of it with the appropriate measures
and decompositions arising from the current C*--context.
First, we transcribe Lemma~\ref{SepChar}
for the current context:
\begin{lem}
\label{CharLf} As before, let $\b f.\in\prod\limits_{n=1}^\infty V_{k_n}$ such that
$\br{\b f.}.\not=0\,.$
Let $\omega$ be a pure state on $\al L.[\b f.],$
and let $\wt{\omega}$ be its strict extension to
the unitaries  $\eta(\R^{(\N)})\subset M\big(\al L.[\b f.]\big)\,.$
Then $\wt{\omega}\circ\eta$ is a character and there exists an
element $\ba \in \R^\N$ with
$\wt{\omega}(\eta(\b x.))=\exp\big(i\bbrk\b x.,{\b a.}.\big)$
for all $\b x.\in\R^{(\N)}$.
\end{lem}
\begin{beweis}
As $\al L.[\b f.]$ is commutative, any pure state $\omega$ of it is a
$*$-homomorphism. Thus the strict extension $\wt{\omega}$ to $\eta(\R^{(\N)})\subset M\big(\al L.[\b f.]\big)$
is also a $*$-homomorphism, hence $\wt{\omega}\circ\eta$ is a character. The restriction of $\wt{\omega}\circ\eta$ to the subgroup
$\R^n\subset\R^{(\N)}$ is still a character, and it is continuous (since it is determined
by the factor $\mathop{\otimes}\limits_{j=1}^nC_0(\R)$ in $\al L.[\b f.]$
which is the group algebra of $\R^n$) hence of the form
$\wt{\omega}\circ\eta(\b x.)=\exp(i\b x.\cdot\b a.^{(n)})$ for some $\b a.^{(n)}\in\R^n\,.$
Since $\wt{\omega}\circ\eta$ is a character on all of $\R^{(\N)},$ the family
${\{\b a.^{(n)}\in\R^n\,\mid\,n\in\N\}}$ is a consistent family, i.e.,
if $n<m$ then $\b a.^{(n)}$ is the first $n$ entries of $\b a.^{(m)}.$
Thus there is an $\b a.\in\R^\N$ such that $\b a.^{(n)}$ is the first $n$ entries of
$\b a.$ for any $n\in\N\,.$ Then
$\wt{\omega}\circ\eta(\b x.)=\exp\big(i\bbrk\b x.,{\b a.}.\big)$
since for any $\b x.\in\R^n\subset\R^{(\N)}$ this restricts to the previous formula
for $\wt{\omega}\circ\eta\,.$
\end{beweis}
Thus there is a
map from the pure states $\wp_P(\al L.[\b f.])$ to $\R^\N$ denoted by
\[
\xi:\wp_P(\al L.[\b f.])\to\R^\N \qquad\hbox{satisfying}\qquad \wt{\varphi}(\eta(\b x.))
=\exp\big(i\bbrk\b x.,\xi(\varphi).\big)
\;\;\forall\;\b x.\in\R^{(\N)}\,,\;\varphi\in\wp_P(\al L.[\b f.])\,.
\]
\begin{teo}
 \label{FullBM}
For each state $\omega$ of $\R^{(\N)}$ there is an
 $\b f.\in\prod\limits_{n=1}^\infty V_{k_n}$ where
$k_n\in\N$ and a unique state $\omega_0\in\wp_0\big(\al L.[\b f.]\big)$
such that  $\eta^*(\omega_0)=\omega\,.$ Then
\begin{itemize}
\item[(i)] there is a regular Borel probability measure $\nu$
on $\wp\big(\al L.[\b f.]\big)$  concentrated on
the pure states $\wp_P\big(\al L.[\b f.]\big)$ such that
\[
 \omega_0(A)=\int_{\wp_P(\al L.[\b f.])}\varphi(A)\,d\nu(\varphi)\quad\forall\,A\in\al L.[\b f.]\,.
\]
\item[(ii)] The probability measure $\wt{\nu}$ on $\R^\N$ given by
$\wt{\nu}:=\xi_*\nu$ is (up to sets of measure zero) the Bochner--Minlos measure for $\omega\,,$ i.e.,
\[
 {\omega}(\b x.)=\int_{\R^\N}\exp\big(i\bbrk\b x.,{\b y.}.\big)\,d\wt{\nu}(\b y.)
\quad\forall\,\b x.\in\R^{(\N)}\,.
\]
\end{itemize}
\end{teo}

\begin{beweis}
Fix an $\omega\in\wp\big(\R^{(\N)}\big).$ Then by
Theorem~\ref{HostLf}(iv) it suffices to show that there is an
 $\b f.\in\prod\limits_{n=1}^\infty V_{k_n}$ such that
$\lim\limits_{k\to\infty}\lim\limits_{n\to\infty}
\left(\Omega_\omega,\,\pi_k^\omega(f_k)\cdots\pi_n^\omega(f_n)\Omega_\omega\right)=1\,.$
However, since there is an approximate identity $\{E_n\}_{n\in\N}$ of $C_0(\R)$
in $\mathop{\cup}\limits_{n=1}^\infty V_n$,
it is possible to choose an $\b f.$ satisfying this limit condition, and we do this as follows.
Since $\lim\limits_{n\to\infty}\pi^\omega_k(E_n)\Omega_\omega=\Omega_\omega,$ choose for each
$n\in\N$ an $f_n := E_{k_n}$ such that
${\big\|\pi^\omega_n(E_{k_n})\Omega_\omega-\Omega_\omega\big\|}\leq 1/ n^2\,.$ Then
for $1<k<n$ we have:
\begin{eqnarray*}
&&\pi^\omega_k(f_k)\cdots\pi_n^\omega(f_n)\Omega_\omega-\Omega_\omega
=\pi^\omega_k(f_k)\cdots\pi_{n-1}^\omega(f_{n-1})\big(\pi_n^\omega(f_n)-\1\big)\Omega_\omega  \\[1mm]
&&\qquad
+\;\pi^\omega_k(f_k)\cdots\pi_{n-2}^\omega(f_{n-2})\big(\pi_{n-1}^\omega(f_{n-1})-\1\big)\Omega_\omega
+\cdots+\big(\pi_k^\omega(f_k)-\1\big)\Omega_\omega\,. \\[1mm]
\hbox{Hence:}\quad &&
\big\|\pi^\omega_k(f_k)\cdots\pi_n^\omega(f_n)\Omega_\omega-\Omega_\omega\big\|
\leq {1\over n^2}+{1\over (n-1)^2}+\cdots+{1\over k^2}  \\[1mm]
&&\hphantom{\big\|\pi^\omega_k(f_k)\cdots\pi_n^\omega(f_n)\Omega_\omega-\Omega_\omega\big\|}
<\int_{k-1}^{n+1}{1\over x^2}\,dx={1\over k-1}-{1\over n+1}
\end{eqnarray*}
from which we see that $\lim\limits_{k\to\infty}\lim\limits_{n\to\infty}
\big\|\pi^\omega_k(f_k)\cdots\pi_n^\omega(f_n)\Omega_\omega-\Omega_\omega\big\|=0,$ and this implies the required
limit condition.

(i)
Since $\al L.[\b f.]$ is separable and commutative, it follows from Theorem~II.2.2 in~\cite{Dav}
that all its GNS--representations are multiplicity free, and hence by Theorem~4.9.4 in~\cite{Ped},
for any state $\omega_0$ on $\al L.[\b f.]$ there is a regular Borel probability measure $\nu$
on $\wp\big(\al L.[\b f.]\big)$  concentrated on
the pure states $\wp_P\big(\al L.[\b f.]\big)$ such that
\[
 \omega_0(A)=\int_{\wp_P(\al L.[\b f.])}\varphi(A)\,d\nu(\varphi)\quad\forall\,A
\in\al L.[\b f.]\,.
\]

(ii)
For the state $\omega_0$ on $\al L.[\b f.],$ let $\wt{\omega}_0$ be its strict extension to
the unitaries $\eta\big(\R^{(\N)}\big)\subset M\big(\al L.[\b f.]\big)\,,$  then  we have
for any countable approximate
identity $\{E_n\}_{n\in\N}\subset\al L._\mu[\b f.]$ that
\begin{eqnarray*}
\wt{\omega}_0\circ\eta(\b x.)&=&\lim_{n\to\infty}\omega_0(\eta(\b x.) E_n)=
\lim_{n\to\infty}\int_{\wp_p(\al L.[\b f.])}\varphi(\eta(\b x.) E_n)\,d\nu(\varphi) \\[1mm]
&=& \int_{\wp_p(\al L.[\b f.])}\lim_{n\to\infty}\varphi(\eta(\b x.)  E_n)\,d\nu(\varphi)
=\int_{\wp_p(\al L.[\b f.])}\wt{\varphi}\circ\eta(\b x.)\,d\nu(\varphi)
\end{eqnarray*}
where we used the Lebesgue Dominated Convergence Theorem in the second line, since
${\big|\varphi(\eta(\b x.) E_n)\big|}\leq1$ and the constant function $1$ is integrable.
If we define a probability measure $\wt{\nu}$ on $\R^\N$ by
$\wt{\nu}:=\xi_*\nu$, where the map
$\xi:\wp_p(\al L.[\b f.])\to\R^\N $ given by
$\wt{\varphi}\circ\eta(\b x.)=\exp\big(i\bbrk\b x.,\xi(\varphi).\big)$ for
$\bx \in\R^{(\N)}$ was mentioned above, we obtain
\[
 \omega(\b x.)=\wt{\omega}_0\circ\eta(\b x.)=\int_{\R^\N}\exp\big(i\bbrk\b x.,{\b y.}.\big)\,d\wt{\nu}(\b y.)
\quad\forall\,\b x.\in\R^{(\N)}.
\]
Hence $\wt{\nu}$ coincides (up to sets of measure zero) with the usual Bochner--Minlos measure on $\R^\N$
by uniqueness of the measure on $\R^\N$ producing this decomposition (cf. Lemma ~7.13.5 in \cite{Bog}).
\end{beweis}

Thus we can interpret the Bochner--Minlos Theorem as
an expression of the pure state space decompositions of the $C^*$--algebras $\al L.[\b f.]\,.$
We will not consider the uniqueness of the measures in the decompositions of the
Bochner--Minlos Theorem, as that is easy to prove.

To understand $\al L.[\b f.]$ at a more concrete level, we consider its spectrum $X\,.$
Since $\al L.[\b f.]$ is commutative, we know $\al L.[\b f.]\cong C_0(X),$
and as each
$\omega\in X$ is a character, we obtain from
Propositions~\ref{repLf} and~\ref{Qchar} that
\[
\omega(L_1\otimes L_2\otimes\cdots)=\lim_{n\to\infty}\omega_1(L_1)\omega_2(L_2)\cdots\omega_n(L_n)q^\ell
\]
for all elementary tensors $L_1\otimes L_2\otimes\cdots\in\br{\b f.^\ell}.\subset\hbox{$*$-alg}\big(\br{\b f.}.\big)$,
where $q\in(0,1]$ and each $\omega_i$ is a character of $\al L.^{\{i\}}=C_0(\R)$ hence a point evaluation
$\omega_i(f)=f(x_i)\,.$
Since $\omega$ is uniquely determined by its values on
$\hbox{$*$-alg}\big(\br{\b f.}.\big)$,
this defines (via Proposition~\ref{Qchar}) a surjective map
\[
 \gamma:\R^{\N}\times(0,1]\to X\cup\{0\}\quad\hbox{by}\quad \gamma(\b x.,q)(L_1\otimes L_2\otimes\cdots)
:=\lim_{n\to\infty}L_1(x_1)L_2(x_2)\cdots L_n(x_n)q^\ell
\]
for $L_1\otimes L_2\otimes\cdots\in\br{\b f.^\ell}.\,.$
To obtain a bijection with $X$ from $\gamma,$ note that
if $A:=L_1\otimes L_2\otimes\cdots=A_1\otimes\cdots\otimes A_m\otimes f_{m+1}^\ell\otimes f_{m+2}^\ell\otimes\cdots
\in\hbox{$*$-alg}\big(\br{\b f.}.\big),$ then
\[
 \prod_{k=1}^{\infty}\,\omega_k(L_k)
=A_1(x_1)\,A_2(x_2)\cdots A_m(x_m)\prod_{k= m+1}^{\infty}\,f_k(x_k)^\ell=0\quad\forall A_i,\,m,\,\ell
\]
if and only if $\lim\limits_{m\to\infty}\prod\limits_{k= m}^{\infty}f_k(x_k)=0\,.$
Thus we define
\[
 N_{\b f.}:=\Big\{\b x.\in\R^{\N}\,\mid\,\lim_{m\to\infty}\prod\limits_{k= m}^{\infty}f_k(x_k)=0\Big\}
\]
and hence the restriction $\gamma:\big(\R^{\N}\backslash N_{\b f.}\big)\times(0,1]\to X$ is a surjection.
That $\gamma$ is bijective, is clear since each $\gamma(\b x.,q)$ is nonzero (as $\b x.\not\in N_{\b f.}$),
and in each factor in the product, a
component of $\al L.[\b f.]$ will separate the characters, and in the last entry, by definition all elementary
tensors will separate different values of $q\,.$ Thus we may identify (as sets) $X$ with
$\big(\R^{\N}\backslash N_{\b f.}\big)\times(0,1]\,.$
Note that $N_{\b f.}$ contains the set $\big\{ \b x.\in\R^{\N}
\,\big|\,x_n\in f_n^{-1}(0)\;\; \hbox{for infinitely many}\;n\big\}\,,$
hence since the $f_n$ are of compact support, $\R^{\N}\backslash N_{\b f.}$ is contained in
the union of sets $\prod\limits_{n=1}^\infty S_n\subset \R^{\N}$ where only finitely many of
the $S_n$ are not relatively compact.

The w$*$-topology of $X$ w.r.t. $\al L.[\b f.]$ is not clear. The most important subset in $X$
is  $X_0:=X\cap{\rm Rep}_0\big(\al L.[\b f.],\C\big)$ which corresponds to
$\big(\R^{\N}\backslash N_{\b f.}\big)\times\{1\}\,.$ We prove that
it is a $G_\delta\hbox{--set.}$ To see this, note that $\omega\in X_0$ if and only if
$\lim\limits_{n\to\infty}\prod\limits_{k=n}^\infty\omega_k(f_k)=1$. This is an increasing limit.
By using approximate identities in each factor $\al L.^{\{k\}}$,
we can find for each $n$ a net
$\{A^{(n)}_\alpha\}\subset\al L.[\b f.],$ $0<A^{(n)}_\alpha<\1\,,$ such that
$\omega\big(\mathop{\overbrace{\1\otimes\cdots\otimes\1}}\limits^{n-1\;{\rm factors}}
\otimes f_n\otimes f_{n+1}\otimes\cdots\big)=\sup\limits_\alpha\omega(A^{(n)}_\alpha)$
for all $\omega\in X\,.$ Define a function $q_{\b f.}:X\to[0,1]$ by
$q_{\b f.}(\omega):=\sup\limits_{\alpha,\,n}\omega(A^{(n)}_\alpha)$ then
$X_0=q_{\b f.}^{-1}(\{1\})\,.$ Since $q_{\b f.}$ is the supremum of continuous functions on $X$
it is lower semicontinous (cf. 6.3 in~\cite{Ko}), i.e.,
$q_{\b f.}^{-1}\big((t,\infty)\big)$ is open for all
$t\in\R\,.$ Since $X_0=q_{\b f.}^{-1}(\{1\})=\mathop{\bigcap}\limits_{n\in\N}q_{\b f.}^{-1}\big((\f n-1,n.,\,\infty)\big),$
it follows that $X_0$ is a $G_\delta\hbox{--set.}$

To make a host algebra out of $\al L.[\b f.],$ i.e., to make $\eta^*$ injective, we need to reduce its
spectrum to $X_0.$ However, since we do not know whether $X_0$ is a locally compact subset of $X$
this is not easy. From the fact that it is a $G_\delta\hbox{--set,}$ we can identify $X_0$ as the common characters of
the decreasing sequence of C*-algebras ${C_0\left(q_{\b f.}^{-1}\big((\f n-1,n.,\,\infty)\big)\right)}\subset\al L.[\b f.],$
where of course $\eta(\R^{(\N)})$ still acts on these as multipliers (i.e., as elements of $C_b(X),$ with pointwise
multiplication).


\section{Hosting the full representation theory of $\R^{(\N)}$}
\label{FullRepRN}

We first want to extend the semi-host algebra $\al L.[\b f.]$ above to an algebra $\al L.\s{\al V.}.$, such that
$\eta^*\left({\rm Rep}\big(\al L.\s{\al V.}.,\al H.\big)\right)={\rm Rep}\big(\R^{(\N)},\al H.\big)\,.$
Recall that for
\[
V_n:=\big\{f\in C_0(\R)\;\big|\;f(\R)\subseteq[0,1],\;\;f\rest[-n,n]=1,\;\;
{\rm supp}(f)\subseteq[-n-1,\,n+1]\big\}\,.
\]
we obtain a multiplicative subsemigroup
$\al V.:=\bigcup\limits_{n=1}^\infty V_n$ in $C_0(\R)$.
Thus, by Theorem~\ref{ProdCls}(iii), $\al V.=\al V.^*,$ implies that
\[
\al A.(\al V.):={\rm Span}\big\{b\in\br{\b f.}.\,\mid\,\b f.\in
\cV^\N\big\}
={\rm Span}\big\{\mathop{\otimes}_{n=1}^\infty g_n\,\mid\,\b g.\sim\b f.\in
\cV^\N\big\}
\]
is a $*$-subalgebra of $\mathop{\bigotimes}\limits_{n=1}^\infty C_0(\R)\,.$
\begin{pro}
 \label{UnivPiL}
There is a $*$-representation $\pi_u:\al A.(\al V.)\to\al B.(\al H._u)$ such that
\[
\pi_u(L_1\otimes L_2\otimes\cdots)=\slim_{n\to\infty}\pi_u^1(L_1)\,\pi_u^2(L_2)\cdots\pi_u^n(L_n)
\]
for all elementary tensors $L_1\otimes L_2\otimes\cdots\in\al A.(\al V.),$ where
$\pi_u^k:C_0(\R)\to\al B.(\al H._u)$ are as before
(cf. text above Definition~\ref{DefLf}).
\end{pro}

\begin{beweis}
By Proposition~\ref{repLf}(iii), $\pi_u$ is already a $*$-representation on each
$\hbox{$*$-alg}\big(\br{\b f.}.\big)$ for $\b f.\in\cV^\N$,
hence it is a linear map on each
$\br{\b f.}.$ for $\b f.\in\cV^\N.$ However, by
Proposition~\ref{prop:2.7}(iv) we know that
for $\b f.,\,\b g.\in\cV^\N$ with $\br{\b f.}.\not=\{0\}\not=\br{\b g.}.$ we have
$\br{\b f.}.\cap\br{\b g.}.=\{0\}$ if and only if $\br{\b f.}.\not\sim\br{\b g.}..$
Thus the set of spaces ${\big\{\br{\b f.}.\,\mid\,\b f.\in
\cV^\N\big\}}$ is labelled by
the equivalence classes ${[\b f.]}\subset\cV^\N,$
and by Proposition~\ref{prop:2.7}(iv), the sum of the subspaces
$\br{\b f.}.$ is direct.
Thus, since $\pi_u$ is defined as a linear map on each $\br{\b f.}.,$ it extends uniquely to a linear map $\pi_u$ on
 $\al A.(\al V.)={\rm Span}\big\{b\in\br{\b f.}.\,\mid\,\b f.\in \cV^\N\big\}.$

To show that $\pi_u$ is a $*$-homomorphism, it suffices to check this on the
elementary tensors $\mathop{\otimes}\limits_{n=1}^\infty g_n$ with
$\bg \sim \bff \in \cV^\N.$
For $\b f.,\,\b g.\in\cV^\N,$
let
\[
 A=A_1\otimes\cdots\otimes A_{k-1}\otimes f_k\otimes f_{k+1}\otimes\cdots\in\br{\b f.}.\quad\hbox{and}\quad
B=B_1\otimes\cdots\otimes B_{k-1}\otimes g_k\otimes g_{k+1}\otimes\cdots\in\br{\b g.}.
\]
where we can choose the same $k$ for both. Then by Proposition~\ref{repLf}(ii) we have
\begin{eqnarray*}
\pi_u(A)&=&\pi_u^1(A_1)\cdots\pi_u^{k-1}(A_{k-1})F_{u,k}[\b f.]\qquad\hbox{and}\qquad
\pi_u(B)=\pi_u^1(B_1)\cdots\pi_u^{k-1}(B_{k-1})\,F_{u,k}[\b g.]  \\[1mm]
\hbox{and}\qquad &&
\pi_u(AB)=\pi_u^1(A_1B_1)\cdots\pi_u^{k-1}(A_{k-1}B_{k-1})F_{u,k}[\b f.\cdot\b g.] \\[1mm]
\hbox{where}\qquad && F_{u,k}[\b f.]:=\slim_{n\to\infty}\pi_u^k(f_k)\cdots\pi_u^n(f_n)\,.
\end{eqnarray*}
Since $\pi_u^j$ is a representation for all $j,$ we only need to show that
$ F_{u,k}[\b f.\cdot\b g.]= F_{u,k}[\b f.]\,F_{u,k}[\b g.]$ to establish that
$\pi_u(AB)=\pi_u(A)\,\pi_u(B)\,.$ We have
\begin{eqnarray*}
 F_{u,k}[\b f.\cdot\b g.]&=&\slim_{n\to\infty}\pi_u^k(f_kg_k)\cdots\pi_u^n(f_ng_n) =\slim_{n\to\infty}\pi_u^k(f_k)\cdots\pi_u^n(f_n)\,\pi_u^k(g_k)\cdots\pi_u^n(g_n)  \\[1mm]
&=&\slim_{n\to\infty}\pi_u^k(f_k)\cdots\pi_u^n(f_n)\cdot\slim_{m\to\infty}\pi_u^k(g_k)\cdots\pi_u^m(g_m)
= F_{u,k}[\b f.]\,F_{u,k}[\b g.]
\end{eqnarray*}
since the operator product is jointly continuous in the strong operator topology on bounded subsets.
Thus $\pi_u$ is a homomorphism. To see that it is a $*$-homomorphism, note that
\[
\pi_u(A)^*=\pi_u^1(A_1^*)\cdots\pi_u^{k-1}(A_{k-1}^*)F_{u,k}[\b f.] =\pi_u(A^*)
\]
since all $\pi_u^j$ are $*$-homomorphisms with commuting ranges, and
$\br{\b f.}.^*=\br{\b f.}^*.=\br{\b f.}..$ Thus $\pi_u$ is a
$*$-homomorphism of $\al A.(\al V.).$
\end{beweis}
As in Section~\ref{ParHRN}, we define
\begin{defi}
\label{DefL}
The $C^*$-algebra $\al L.\s{\al V.}.$ is the $C^*$-completion of $\pi_u\big(\al A.(\al V.)\big)$ in $\al B.(\al H._u)\,.$
\end{defi}
Note that $\al L.\s{\al V.}.={\rm C}^*\big\{\al L.[\b f.]\,\mid\,\b f.\in\cV^\N\big\}
\subset\al B.(\al H._u)\,.$

We extend the unitary embeddings $\eta:\R^{(\N)}\to UM\big(\al L.[\b f.]\big)$
from above to $\al L.\s{\al V.}.$ as follows.
Define $\eta:\R^{(\N)}\to M\big(\al L.\s{\al V.}.\big)$, where
\begin{eqnarray*}
&&\!\!\!\!\!\!\!
 \eta(x_1,\ldots,x_n,0,0,\ldots)(L_1\otimes L_2\otimes\cdots)=\eta_1(x_1)L_1\otimes\cdots\otimes
\eta_n(x_n)L_n\otimes L_{n+1}\otimes L_{n+2}\otimes\cdots \\[1mm]
&&=\zeta\s{\{1,\ldots,n\}}.(x_1,\ldots,x_n)\big(L_1\otimes  L_2
\otimes \cdots\big)
\end{eqnarray*}
for all $(x_1,\ldots,x_n)\in\R^n,$ $n\in\N,$ $L_i\in\al L.^{\{i\}}=C_0(\R),$
and where $\eta_i:\R\to M({\rm C}^*(\R))$ is the usual unitary embedding.
Clearly, $\eta$ restricts to the previous definition of it on each
$\al L.[\b f.]\subset\al L.\s{\al V.}..$
Then the map $\eta^*:{\rm Rep}\big(\al L.\s{\al V.}.,\al H.\big)\to
{\rm Rep}\big(\R^{(\N)},\al H.\big)$ consists of the strict extension of
(non-degenerate) representations of $\al L.\s{\al V.}.$ to $\eta(\R^{(\N)}),$ i.e.,
\[
 (\eta^*\pi)(\b x.):=\slim_{\alpha\to\infty}\pi\big(\eta(\b x.)E_\alpha\big)\quad
\forall\;\b x.\in\R^{(\N)}\quad\hbox{and any approximate identity}\quad
\{E_\alpha\}\s\alpha\in\Lambda.\subset\al L.\s{\al V.}.\,,
\]
and $\eta^*$ obviously takes irreducibles to irreducibles by commutativity.
\begin{defi}
\label{Rep0L}
Let ${\rm Rep}_0\big(\al L.\s{\al V.}.,\al H.\big)$ denote those non-degenerate $*$-representations
$\pi:\al L.\s{\al V.}.\to\al B.(\al H.)$ for which  $\pi\restriction\al L.[\b f.]
\in{\rm Rep}_0\big(\al L.[\b f.],\al H._{\b f.}\big)$ for all $\b f.,$ where
$\al H._{\b f.}:=\pi\big(\al L.[\b f.]\big)\al H.\,.$ That is, each restriction of $\pi$ to
$\al L.[\b f.]$ has excess operator $Q_\bff=\1$ on its essential subspace $\al H._{\b f.}.$
\end{defi}
 By Proposition~\ref{repLf}, this means that all
\[
 Q_{\b f.}(\pi):=\slim_{n\to\infty}B_n[\b f.]\quad\hbox{are projections, where:}\quad
B_n[\b f.]:=\wt\pi\big(\mathop{\overbrace{\1\otimes\cdots\otimes\1}}\limits^{n-1\;{\rm factors}}
\otimes f_n\otimes f_{n+1}\otimes\cdots\big).
\]
In fact, the projections $Q_{\b f.}(\pi)$ must be the range projections
$P_\pi[\b f.]=\slim\limits_{k\to\infty}F_{\pi,k}^{(1)}$ where
$F_{\pi,k}^{(1)}:=\slim\limits_{n\to\infty}\pi_k(f_k)\cdots\pi_n(f_n)\,.$
Note that a direct sum of representations $\pi_i\in{\rm Rep}_0\big(\al L.\s{\al V.}.,\al H._i\big),$ $i\in I$ (an index set)
is again of the same type, i.e., $\bigoplus\limits_{i\in I}\pi_i\in{\rm Rep}_0\big(\al L.\s{\al V.}.,\bigoplus\limits_{i\in I}\al H._i\big).$
\begin{teo}
 \label{HostL}
Given the preceding notation, we have that
\begin{itemize}
 \item[(i)] $\eta:\R^{(\N)}\to M\big(\al L.\s{\al V.}.\big)$ is continuous w.r.t. the strict topology of
$M\big(\al L.\s{\al V.}.\big)\,.$
\item[(ii)]
The map $\eta^*$ is injective on ${\rm Rep}_0\big(\al L.\s{\al V.}.,\al H.\big)\,.$
\item[(iii)] The range $\eta^*\left({\rm Rep}\big(\al L.\s{\al V.}.,\al H.\big)\right)$
is the same as $\eta^*\left({\rm Rep}_0\big(\al L.\s{\al V.}.,\al H.\big)\right)$
and is all of ${\rm Rep}\big(\R^{(\N)},\al H.\big)\,.$
\item[(iv)] $\pi\in{\rm Rep}_0\big(\al L.\s{\al V.}.,\al H.\big)$ if and only if $\pi$ is normal w.r.t. $\pi_u.$
\end{itemize}
\end{teo}
\begin{beweis}
(i) Since $\eta:\R^{(\N)}\to M\big(\al L.\s{\al V.}.\big)$ is bounded, it suffices to show that the space
\[
\big\{L\in\al L.\s{\al V.}.\,\mid \hbox{the map}\quad\R^{(\N)} \ni \bx \mapsto \eta(\bx)L\in\cL\s{\al V.}.
\quad\hbox{is norm continuous}\big\}
\]
is dense in $\al L.\s{\al V.}..$ But this follows from the fact that by Theorem~\ref{HostLf}(i), this space contains all $\br{\b f.}.\subset\al L.[\b f.],$ and these spaces span $\al A.(\al V.)$ which is dense in
$\al L.\s{\al V.}..$

(ii)
 Consider
 $\pi,\,\pi'\in{\rm Rep}_0\big(\al L.\s{\al V.}.,\al H.\big)$ such that $\eta^*\pi=\eta^*\pi'.$
Then for the restrictions to $\R^n\subset\R^{(\N)}$ we have
$\wt\pi^{(n)}:=\eta^*\pi\restriction\R^n=\eta^*\pi'\restriction\R^n:={\wt{\pi}'}{}^{(n)}\,.$
Moreover, $\al L.^{(n)}$ embeds in $M(\al L.\s{\al V.}.)$ as $\al L.^{(n)}\otimes\1$ (acting on the elementary tensors),
hence $\pi$ also extends to it to define a non-degenerate $\pi^{(n)}:\al L.^{(n)}\to \al B.(\al H._\pi)\,.$
Since $\eta$ is defined via the natural actions, we have
$\eta(\b x.)\al L.^{(n)}\subseteq\al L.^{(n)}$ for all $\b x.\in\R^n\,.$
Since
\[
\wt\pi^{(n)}(\b x.)\, \pi^{(n)}(L)\,\pi(A)=\eta^*\pi(\b x.)\,\pi(LA)=\pi\big(\eta(\b x.)LA\big)=
\pi^{(n)}(\eta(\b x.)L)\,\pi(A)
\]
for all $\b x.\in\R^n,\; L\in\al L.^{(n)},\;A\in\al L.\s{\al V.}.,$ we see by nondegeneracy of $\pi$ that
$\wt\pi^{(n)}(\b x.)\, \pi^{(n)}(L)=\pi^{(n)}(\eta(\b x.)L)$ for all $L\in\al L.^{(n)},$
and hence since $\wt\pi^{(n)}$ and
$\pi^{(n)}$ are non-degenerate and $\al L.^{(n)}$ is a host algebra for $\R^n,$ this relation
gives a bijection between $\wt\pi^{(n)}$ and $\pi^{(n)}.$ We conclude from
$\wt\pi^{(n)}=  {\wt{\pi}'}{}^{(n)}$ that
$\pi^{(n)}={\pi'}^{(n)}$ for all $n\,.$
A similar argument for the $k^{\rm th}$component alone, also shows that $\pi_k=\pi'_k$ for all $k\,.$
Now for each elementary tensor
$L_1\otimes L_2\otimes\cdots\in\hbox{$*$-alg}\big(\br{\b f.}.\big)\subset\al L.[\b f.]$
we know by Proposition~\ref{repLf}(iii) that
\begin{equation}
\label{PiQ}
 \pi(L_1\otimes L_2\otimes\cdots)=\slim_{n\to\infty}\pi^{(n)}\big(L_1\otimes L_2\otimes\cdots\otimes L_n\big)
Q_{\b f.}(\pi)\,.
\end{equation}
Recall that by hypothesis  $Q_{\b f.}(\pi)=P_\pi[\b f.]
=\slim\limits_{k\to\infty}F_{\pi,k}^{(1)}$, where
$F_{\pi,k}^{(1)}:=\slim\limits_{n\to\infty}\pi_k(f_k)\cdots\pi_n(f_n)\,.$
Analogous expressions hold for $\pi',$ thus
 since $\pi_k=\pi'_k$ for all $k\,,$ it follows that
$Q_{\b f.}(\pi)=Q_{\b f.}(\pi')$ and hence from Equation~(\ref{PiQ})
it follows from $\pi^{(n)}={\pi'}^{(n)}$ for all $n\,,$ that
 $\pi$ and $\pi'$
coincides on all $\al L.[\b f.]$ hence on all of $\al L.\s{\al V.}.,$ which proves the claim.

(iii) Let $\pi\in{\rm Rep}_0\big(\al L.\s{\al V.}.,\al H.\big)$ and let $\wt\pi$ be its strict extension
to $M\big(\al L.\s{\al V.}.\big).$ As $\wt\pi$ is strictly continuous, (i) implies that the
unitary representation $\eta^*(\pi)=\wt{\pi}\circ\eta:\R^{(\N)}\to\al U.(\al H.)$ is strong operator
continuous, i.e., $\eta^*\left({\rm Rep}\big(\al L.\s{\al V.}.,\al H.\big)\right)\subseteq {\rm Rep}\big(\R^{(\N)},\al H.\big)\,.$
To prove the claim of this theorem, we need to show that for each $\pi\in{\rm Rep}\big(\R^{(\N)},\al H._\pi\big),$
there is a $\pi_{(0)}\in{\rm Rep}_0\big(\al L.\s{\al V.}.,\al H._\pi\big)$ such that
$\eta^*\pi_{(0)}=\pi\,.$ Since each $\pi\in{\rm Rep}\big(\R^{(\N)},\al H._\pi\big)$ is a
direct sum of cyclic representations, and $\eta^*$ preserves direct sums, it suffices to
show that for each cyclic $\pi\in{\rm Rep}\big(\R^{(\N)},\al H._\pi\big),$
there is a $\pi_{(0)}\in{\rm Rep}_0\big(\al L.\s{\al V.}.,\al H._\pi\big)$ such that
$\eta^*\pi_{(0)}=\pi\,.$
Fix a cyclic $\pi\in{\rm Rep}\big(\R^{(\N)},\al H._\pi\big),$ then there is a
projection $P_\pi\in \pi_u(\R^{(\N)})'$ such that
$\pi=(P_\pi\pi_u)\restriction\al H._\pi$ where $\al H._\pi=P_\pi\al H._u\,.$
Recall the inclusion
$\R\to \R^{(\N)}, x \mapsto x e_k$,
so let $\pi_k:\R\to \al U.(\al H._\pi)$ be $\pi_k(x):=  \pi(xe_k)$.
 By the host algebra property of $C^*(\R)\cong C_0(\R),$ this
produces a unique non-degenerate representation $\pi_k:C_0(\R)\to \al B.(\al H._u)\,,$ which is
characterized by $\pi_k(x)\,\pi_k(L)=\pi_k(\eta_k(x)L)=\pi\big(\eta(0,\ldots,0,x,0,0,\ldots)(\1\otimes
\cdots\1\otimes L\otimes\1\otimes\cdots)\big)$ ($x$ and $L$ in the $k^{\rm th}$entries) for all
$x\in\R$ and $L\in C_0(\R).$ Since
\begin{eqnarray*}
&&\pi\big(\eta(0,\ldots,0,x,0,\ldots)(\1\otimes\cdots\1\otimes L\otimes\1\otimes\cdots)\big) \\[1mm]
&&\qquad= P_\pi \pi_u\big(\eta(0,\ldots,0,x,0,\ldots)(\1\otimes\cdots\1\otimes L\otimes\1\otimes\cdots)\big)
\restriction\al H._\pi \\[1mm]
&&\qquad= P_\pi \pi_u^k(x)\,\pi_u^k(L)\restriction\al H._\pi
=\pi_k(x)\,P_\pi\pi_u^k(L)\restriction\al H._\pi
\end{eqnarray*}
we get that $\pi_k(L)=P_\pi\pi_u^k(L)\restriction\al H._\pi$ for all $L\in C_0(\R).$
Since the set of representations
${\big\{\pi_k:C_0(\R)\to \al B.(\al H._\pi)\,\mid\,k\in\N\big\}}$ have commuting ranges,
we can apply Lemma~\ref{RepSet} (with the choice $Q=\1$) to define a representation
$\pi_{(0)}:\hbox{$*$-alg}\big(\br{\b f.}.\big)\to\al B.(\al H._\pi),$ for all $\b f.,$ and we need to show that
$\pi_{(0)}$ extends to a representation of  $\al L.\s{\al V.}..$
Now $P_\pi$ commutes with the images of all $\pi_u^k$ (since it commutes with $\pi_u(\R^{(\N)})$) hence
all $\pi_u^k(L)$ preserve $\al H._\pi$  and so by its definition $\pi_u(\al L.\s{\al V.}.)$ preserves $\al H._\pi.$
Thus the map $A\in\al L.\s{\al V.}.\to P_\pi\pi_u(A)\restriction\al H._\pi$ is a $*$-representation of $\al L.\s{\al V.}.$
and it coincides with $\pi_{(0)}$ on each
 $\hbox{$*$-alg}\big(\br{\b f.}.\big)$ because
\begin{eqnarray*}
 P_\pi\pi_u(L_1\otimes L_2\otimes\cdots)\restriction\al H._\pi&=&
\slim_{n\to\infty}P_\pi\pi_u^1(L_1)\,\pi_u^2(L_2)\cdots\pi_u^n(L_n)\restriction\al H._\pi \\[1mm]
&=&\slim_{n\to\infty}\pi_1(L_1)\,\pi_2(L_2)\cdots\pi_n(L_n)=\pi_{(0)}(L_1\otimes L_2\otimes\cdots)
\end{eqnarray*}
for all elementary tensors $L_1\otimes L_2\otimes\cdots\in\al A.(\al V.).$
This defines a $*$-representation
\break $\pi_{(0)}:\al L.\s{\al V.}.\to\al B.(\al H._\pi)$                
by $\pi_{(0)}(A)=P_\pi\pi_u(A)\restriction\al H._\pi$ for all $A\in\al L.\s{\al V.}..$
To see that it is non-degenerate, note that its restriction to any
$\al L.[\b f.]\subset\al L.\s{\al V.}.$ has essential projection
$P_\pi[\b f.]=\slim\limits_{k\to\infty}\wt{F}_k$ where
$\wt{F}_k:=\slim\limits_{n\to\infty}\pi_k(f_k)\cdots\pi_n(f_n)$
by Theorem~\ref{HostLf}(iii) and Lemma~\ref{RepSet}(ii).
It is suffices to show that for each nonzero $\psi\in\al H._\pi$
there is a sequence $\b f.\in\cV^\N$
such that $P_\pi[\b f.]\psi\not=0.$ Fix a nonzero $\psi\in\al H._\pi.$
Since there is an approximate identity of $C_0(\R)$ in $\al V.,$ it
is possible to choose for each $n\in\N$ a $f_n\in\al V.$ such that
${\|\psi-\pi_n(f_n)\psi\|}<1/n^2,$ hence we may write
$\pi_n(f_n)\psi=\psi +\xi_n/n^2$ where $\|\xi_n\|\leq 1.$ Then
\begin{eqnarray*}
\pi_k(f_k)\cdots\pi_n(f_n)\psi&=& \psi +{1\over n^2}\,\pi_k(f_k)\cdots\pi_{n-1}(f_{n-1})\,\xi_n  \\[1mm]
&&\quad + \;{ 1\over(n-1)^2}\,\pi_k(f_k)\cdots\pi_{n-2}(f_{n-2})\,\xi_{n-1}+\cdots
+{ 1\over k^2}\,\xi_k\;.
\end{eqnarray*}
Thus
\[ \wt{F}_k\psi= \psi+\sum_{j=k}^\infty{1\over j^2}\,\prod_{\ell=k}^{j-1}\pi_\ell(f_\ell)\,\xi_j,
\qquad\hbox{where}\qquad \big\|\prod_{\ell=k}^{j-1}\pi_\ell(f_\ell)\,\xi_j\big\|\leq 1 \]
and hence $P_\pi[\b f.]\psi=\psi$ as the series converges.
Thus $\pi_{(0)}$ is non-degenerate.

Since
$\pi_{(0)}$ is obviously normal to $\pi_u,$ it follows that
the excess operator is $Q=\1$ for the restriction of $\pi_{(0)}$ to any
$\al L.[\b f.],$ and hence $\pi_{(0)}\in{\rm Rep}_0\big(\al L.\s{\al V.}.,\al H._\pi\big)\,.$
To see that $\eta^*\pi_{(0)}=\pi\,,$ note that for $\b x.\in\R^k\subset\R^{(\N)}$ we have
\begin{eqnarray*}
&&\eta^*\pi_{(0)}(\b x.)\,\pi_{(0)}(L_1\otimes L_2\otimes\cdots)=
\pi_{(0)}\big(\eta(\b x.)(L_1\otimes L_2\otimes\cdots)\big)\\[1mm]
&&\qquad=\pi_{(0)}\big(\eta_1(x_1)L_1\otimes\cdots\otimes\eta_k(x_k)L_k\otimes L_{k+1}\otimes L_{k+2}\otimes\cdots\big)\\[1mm]
&&\qquad= P_\pi\pi_u\big(\eta_1(x_1)L_1\otimes\cdots\otimes\eta_k(x_k)L_k\otimes L_{k+1}\otimes
L_{k+2}\otimes\cdots\big)\restriction\al H._\pi\\[1mm]
&&\qquad= P_\pi\,\slim_{n\to\infty}\pi_u^1(\eta_1(x_1)L_1)\,\pi_u^2(\eta_2(x_2)L_2)\cdots\pi_u^k(\eta_k(x_k)L_k)
\,\pi_u^{k+1}(L_{k+1})\cdots\pi_u^n(L_n)\restriction\al H._\pi \\[1mm]
&&\qquad= P_\pi\,\slim_{n\to\infty}\pi_u^1(x_1)\pi_u^1(L_1)\,\pi_u^2(x_2)\pi_u^2(L_2)\cdots\pi_u^k(x_k)\pi_u^k(L_k)
\,\pi_u^{k+1}(L_{k+1})\cdots\pi_u^n(L_n)\restriction\al H._\pi \\[1mm]
&&\qquad=\pi_1(x_1)\cdots\pi_k(x_k) \,\slim_{n\to\infty}\pi_1(L_1)\,\pi_2(L_2)\cdots
\pi_n(L_n) \\[1mm]
&&\qquad=\pi(\b x.)\pi_0(L_1\otimes L_2\otimes\cdots)
\end{eqnarray*}
for all elementary tensors $L_1\otimes L_2\otimes\cdots\in\al A.(\al V.),$
hence $\eta^*\pi_0(\b x.)\,\pi_{(0)}(A)=\pi(\b x.)\,\pi_{(0)}(A)$ for all $A\in\al L.\s{\al V.}..$
Since $\pi_{(0)}$ is non-degenerate, it follows that $\eta^*\pi_{(0)}=\pi$ as required.

(iv) By Theorem~\ref{HostLf}(v) we have that $\pi\in{\rm Rep}_0\big(\al L.\s{\al V.}.,\al H.\big)$ if and only if $\pi\restriction
\al L.[\b f.]$ (on its essential subspace $\al H._{\b f.})$ is normal w.r.t.
$\pi_u\big(\al L.[\b f.]\big)$ for all $\b f.\in\cV^\N.$
Let $\pi\in{\rm Rep}\big(\al L.\s{\al V.}.,\al H.\big)$ be
normal  w.r.t. $\pi_u(\al L.\s{\al V.}.)$. Then it is continuous
on bounded sets w.r.t. the strong operator topologies
of both sides, and it follows that this is true for its restrictions to each $\pi_u(\al L.[\b f.]),$
and hence that each restriction is normal w.r.t. $\pi_u.$ Thus $\pi\in{\rm Rep}_0\big(\al L.\s{\al V.}.,\al H.\big).$

Conversely, given $\pi\in{\rm Rep}_0\big(\al L.\s{\al V.}.,\al H.\big)$ then
by part (iii) $\eta^*\pi$ is a continuous representation
of $\R^{(\N)},$ and by
Proposition~\ref{repLf}(iii) (with $Q=\1)$ we have that on each $\al H._{\b f.}$
\[
\pi(L_1\otimes L_2\otimes\cdots)=\slim_{n\to\infty}\pi^{(n)}\big(L_1\otimes L_2\otimes\cdots\otimes L_n\big)
=\slim_{n\to\infty}\pi_1(L_1)\pi_2(L_2)\cdots\pi_n( L_n)
\]
for all elementary tensors $L_1\otimes L_2\otimes\cdots\in\br{\b f.^\ell}.\subset\hbox{$*$-alg}\big(\br{\b f.}.\big).$
 Now $\pi_u(\R^{(\N)})''={\{\pi_u^{(n)}(\al L.^{(n)})\,\mid\,n\in\N\}''}=
\pi_u\left(\big\{\al L.[\b f.]\,\mid\,\b f.\in
\cV^\N\big\}\right)''=\pi_u(\al L.\s{\al V.}.)''$ and a similar equation holds for $\pi.$
Since the cyclic components of $\pi$ are contained in the direct summands of $\pi_u,$
there is a normal map $\varphi:\pi_u(\al L.\s{\al V.}.)''\to\al B.(\al H.)$ such that
$\varphi\circ\pi_u=\pi.$ Thus $\pi$ is normal to $\pi_u.$
\end{beweis}
Thus $\al L.\s{\al V.}.$ is a
semi-host for the full representation theory of $\R^{(\N)},$
i.e., $\eta^*:{\rm Rep}\big(\al L.\s{\al V.}.,\al H.\big)\to{\rm Rep}\big(\R^{(\N)},\al H.\big)$
is surjective, but not necessarily injective.
We want to examine the remaining representations of $\al L.\s{\al V.}.$ outside of
${\rm Rep}_0\big(\al L.\s{\al V.}.,\al H.\big)\,.$ Denote the universal representation
of $\al L.\s{\al V.}.$ by $\pi\s U.:\al L.\s{\al V.}.\to\al B.(\al H._U)$ (not to be confused with
the defining representation $\pi_u)\,.$ Let
\[
 \al Q.:=\big\{Q_{\b f.}(\pi\s U.)\,\mid\,\b f.\in \cV^\N\big\}
\subset\al L.\s{\al V.}.'':=\pi\s U.(\al L.\s{\al V.}.)''\,,
\]
i.e., the set of all excess operators w.r.t. $\pi\s U.\,.$
Since $\al Q.$ is in the positive part of the unit ball of $\al L.\s{\al V.}.''$,
it has a natural partial order, and in a moment we will see that $\al Q.$ is a
multiplicative semigroup. Let
\[
{\rm Rep}\big(\al Q.,\al H.\big):=\big\{\gamma:\al Q.\to\al B.(\al H.)\,\mid\,
0\leq\gamma(Q_1)\leq\1,\;\;\gamma(Q_1Q_2)=\gamma(Q_1)\gamma(Q_2)\;\forall\, Q_i\in\al Q.\big\}\,.
\]
\begin{pro}
\label{repL}
  With notation above, we have
\begin{itemize}
 \item[(i)] $Q_{\b f._1}(\pi\s U.)\cdot Q_{\b f._2}(\pi\s U.)
=Q_{\b f._1\cdot\b f._2}(\pi\s U.)$ for all $\b f._i\in
\cV^\N\,.$
Thus $\al Q.$ is a multiplicative semigroup, and the map $\eqc{\b f.}.\to
Q_{\b f.}(\pi\s U.)$ is a surjective homomorphism $\al V._\infty\to\al Q.$ of multiplicative
semigroups where $\al V._\infty:=\{\eqc{\b f.}.\,\mid\,\b f.\in
\cV^\N\}\,.$
\item[(ii)] Fix a non-degenerate $*$-representation
$\pi:\al L.\s{\al V.}.\to\al B.(\al H._\pi)$. Then the map $\eqc{\b f.}.\to Q_{\b f.}(\pi)$
defines a representation of $\al V._\infty$
as well as of $\al Q.\,.$
Thus every $\pi\in{\rm Rep}\big(\al L.\s{\al V.}.,\al H.\big)$ is of the form:
\[
 \pi(A):=\pi_0(A)\gamma(\b f.)\quad\hbox{for}\quad
A\in \br{\b f.}.
\]
for some $\pi_0\in{\rm Rep}_0\big(\al L.\s{\al V.}.,\al H.\big)$ and
$\gamma\in{\rm Rep}\big(\al Q.,\al H.\big)$ with $\gamma(\al Q.)\subset\pi(\al L.\s{\al V.}.)'\,.$
\end{itemize}
\end{pro}
\begin{beweis}
(i) Recall that  $Q_{\b f.}(\pi):=\slim\limits_{n\to\infty}B_n[\b f.]$, where
$B_n[\b f.]:=\wt\pi\big(\mathop{\overbrace{\1\otimes\cdots\otimes\1}}\limits^{n-1\;{\rm factors}}
\otimes f_n\otimes f_{n+1}\otimes\cdots\big).$
Since the operator product is jointly continuous on bounded subsets we have:
\begin{eqnarray*}
Q_{\b f.}(\pi\s U.)\cdot Q_{\b g.}(\pi\s U.)&=&
\slim_{n\to\infty} \wt{\pi\s U.}\big(\1\otimes\cdots\otimes\1\otimes
f_n\otimes f_{n+1}\otimes\cdots\big)
\slim_{k\to\infty} \wt{\pi\s U.}\big(\1\otimes\cdots\otimes\1\otimes g_k\otimes g_{k+1}\otimes\cdots\big)
\\[1mm]
&=&\slim_{n\to\infty} \wt{\pi\s U.}\big(\1\otimes\cdots\otimes\1\otimes f_n\otimes f_{n+1}\otimes\cdots\big)
\, \wt{\pi\s U.}\big(\1\otimes\cdots\otimes\1\otimes g_n\otimes g_{n+1}\otimes\cdots\big)
\\[1mm]
&=&\slim_{n\to\infty} \wt{\pi\s U.}\big(\1\otimes\cdots\otimes\1\otimes f_ng_n\otimes f_{n+1}g_{n+1}\otimes\cdots\big)\\[1mm]
&=& Q_{\b f.\cdot\b g.}(\pi\s U.)\,.
\end{eqnarray*}
It will suffice for this part to show that the map $[\b f.]\to Q_{\b f.}(\pi\s U.)$ is well-defined,
i.e., that $Q_{\b f.}(\pi\s U.)$ only depends on the equivalence class $[\b f.]$ not on any particular representative
which is chosen. However, this is immediate from the definition of $Q_{\b f.}(\pi\s U.)\,.$

(ii) By the universal property of $\pi\s U.$ (cf. Theorem~10.1.12 in \cite{KR2}) there is a central projection
$P_\pi\in Z\big(\pi\s U.(\al L.\s{\al V.}.)''\big)$ and a *-isomorphism of Von Neumann algebras
$\alpha:P_\pi\pi\s U.(\al L.\s{\al V.}.)''\to \pi(\al L.\s{\al V.}.)''$  such that
$\pi(A)=\alpha(P_\pi\pi\s U.(A))$ for all $A\in\al L.\s{\al V.}..$ The map $\alpha$ is normal in both directions
(cf. Proposition~2.5.2 in \cite{Ped}). It is also true that
$\wt\pi(A)=\alpha(P_\pi\wt\pi\s U.(A))$ for all $A\in M(\al L.\s{\al V.}.).$ So it follows from
\begin{eqnarray*}
Q_{\b f.}(\pi)&=&
\slim_{n\to\infty} \wt{\pi}\big(\1\otimes\cdots\otimes\1\otimes
f_n\otimes f_{n+1}\otimes\cdots\big)=\alpha\left(P_\pi\slim_{n\to\infty} \wt{\pi\s U.}\big(\1\otimes\cdots\otimes\1\otimes
f_n\otimes f_{n+1}\otimes\cdots\big)\right)   \\[1mm]
&=&\alpha(P_\pi Q_{\b f.}(\pi\s U.))
\end{eqnarray*}
and part (i) that $Q_{\b f.}(\pi)\cdot Q_{\b g.}(\pi)=Q_{\b f.\cdot\b g.}(\pi)$ hence
the map $\eqc{\b f.}.\to Q_{\b f.}(\pi)$
defines a representation of $\al V._\infty$
as well as of $\al Q.\,.$ The second claim is immediate.
%
\end{beweis}
Thus the  the additional part of
${\rm Rep}\big(\al L.\s{\al V.}.,\al H.\big)$ to ${\rm Rep}\big(\R^{(\N)},\al H.\big)$
is in ${\rm Rep}\big(\al Q.,\al H.\big)\,.$

By definition, each $Q\in\al Q.$ is the strong operator limit of increasing positive elements
in $\pi_U(\al L.\s{\al V.}.),$ so it is a lower semi-continuous function on the spectrum of
$\al L.\s{\al V.}..$ In fact, $\al Q.$ is in the monotone closure
$\al L.\s{\al V.}.^m$
(cf. \cite[Thm.~6.8 and above, p.~182]{Tak}).
Let $X$ be the spectrum of $\al L.\s{\al V.}.,$ and let $X_0:=X\cap{\rm Rep}_0\big(\al L.\s{\al V.}.,\C\big).$
Then since $\omega(Q)$ must be idempotent for $\omega\in X_0,$ $Q\in\al Q.,$ it has to be $0$ or $1.$
Thus $X_0\subset Q^{-1}(\{0\})\cup Q^{-1}(\{1\}),$ and by the definition of ${\rm Rep}_0\big(\al L.\s{\al V.}.,\C\big)$
we get that
\[
X_0=\bigcap_{Q\in\al Q.}\left(Q^{-1}(\{0\})\cup Q^{-1}(\{1\})\right)\,.
\]
This suggests that to obtain a full host algebra for $\R^{(\N)}$ we only need to apply the homomorphism
which factors out by $\bigcup\limits_{Q\in\al Q.}Q^{-1}((0,1))\,,$ but this is not possible, because we do not
know whether the last set is open, as the $Q$ are only lower semi-continuous.

\section{Discussion}

Here we constructed an infinite tensor product of the algebras $C_0(\R),$
denoted $\al L.\s{\al V.}.,$ and used it to
find semi--hosts for the full continuous representation theory of $\R^{(\N)}.$
Due to commutativity, these were as useful as host algebras, because $\eta^*$
preserves irreducibility in this context. We also
interpreted the Bochner--Minlos theorem in $\R^{(\N)}$ as the pure state space decomposition
of the partial hosts which $\al L.\s{\al V.}.$ comprises of. We analyzed the representation theory
of $\al L.\s{\al V.}.,$  and showed that $\eta^*$ is a bijection between
${\rm Rep}_0\big(\al L.\s{\al V.}.,\al H.\big)$ and
 ${\rm Rep}\big(\R^{(\N)},\al H.\big)\,,$ but that there is an extra part which
 essentially consists of the representation theory of a multiplicative semigroup
 $\al Q..$

 Much further analysis remains, e.g. the topological structure of the
 spectrum $X$ of $\al L.\s{\al V.}.,$ especially the important question as to whether
 $X_0$ is locally compact with the relative topology. Moreover, one can easily apply
 the methods developed here to construct infinite tensor products of other C*-algebras
 without nontrivial projections. A very important issue, is to extend the C*-algebraic
 interpretation of the Bochner--Minlos theorem developed here, to general nuclear spaces.

\section*{Acknowledgements}

HG wishes to thank the Mathematics Department of the
Technische Universit\"at Darmstadt
who generously supported his visit to Darmstadt in 2009,
as well as the Deutscher Akademischer Austauschdienst (DAAD) who funded his trip to Germany.

\bigskip

\providecommand{\bysame}{\leavevmode\hbox to3em{\hrulefill}\thinspace}

\end{document}